\documentclass[final,1p,times,11pt]{article}

\usepackage{enumerate}
\usepackage{epsfig,alltt}
\usepackage[english]{babel}
\usepackage{blindtext}
\usepackage{dsfont}
\usepackage{amssymb}
\usepackage{amsfonts}
\usepackage{graphicx}
\usepackage{amsmath}
\usepackage{amsthm}
\usepackage{natbib}
\usepackage[utf8]{inputenc}

\usepackage[left=2.5cm,right=2.5cm,top=2.5cm,bottom=2.5cm]{geometry}

\usepackage{floatflt}
\usepackage{makeidx}
\usepackage{algorithm}

\usepackage{algorithmic}

\usepackage{url}
\usepackage{float} 
\theoremstyle{remark}

\usepackage{color}

\newcommand{\begeq}[1]{\begin{equation} \label{#1}}
\newcommand{\fineq}{\end{equation}}

\title{Recursive kernel density estimators under missing data}
\author{Yousri Slaoui\\
Universit\'e de Poitiers, Laboratoire de Math\'ematiques et Application,\\ Futuroscope Chasseneuil, France\\
E.mail: Yousri.Slaoui@math.univ-poitiers.fr}

\begin{document}

\newtheorem{theor}{Theorem}
\newtheorem{prop}{Proposition}
\newtheorem{lemma}{Lemma}
\newtheorem{lem}{Lemma}
\newtheorem{coro}{Corollary}

\newtheorem{prof}{Proof}
\newtheorem{defi}{Definition}
\newtheorem{rem}{Remark}

\date{}
\maketitle

\textit{Abstract}: In this paper we propose an automatic bandwidth selection of the recursive kernel density estimators with missing data in the context of global and local density estimation. We showed that, using the selected bandwidth and a special stepsize, the proposed recursive estimators outperformed the nonrecursive one in terms of estimation error in the case of global estimation. However, the recursive estimators are much better in terms of computational costs. We corroborated these theoretical results through simulation studies and on the simulated data of the Aquitaine cohort of $\texttt{HIV}$-1 infected patients and on the coriell cell lines using the chromosome number 11.\\

2010 \textit{Mathematics Subject Classification}: Primary 62G07, 62L20; Secondary 65D10, 62N01.\\
\textit{Key words and phrases:} Density estimation; Stochastic approximation algorithm; Smoothing, curve fitting; missing data.

\section{Introduction}
Kernel density estimator (KDE) were introduced in~\citet{Par62} and~\citet{Ros56} have been widely used in many applications research including clinical trials, epidemiology, genetics as an exploratory tool.\\
There has been an intensive work in the literature about the bandwidth selection methods of this estimator. See for example~\citet{Dui76},~\citet{Rud82},~\citet{Sco87} and~\citet{Mar88}.
However, in many situation the full data are not available and are subjects to missing data. Several methods have been used to overcome the problem and can be found in the literature. For example, the EM algorithm introduced in~\citet{Dem77}, can be used to obtain maximum likelihood estimates when some data are missing, and multiple imputation, which is a Monte Carlo technique in which the missing values are replaced by simulated versions, and the results are combined to produce estimates which incorporate missing data, see the paper of~\citet{Lit02}. Note that all of these methods are concerned with the estimation of unknown parameters.\\
In the framework of nonparametric estimation,~\citet{Dub09} used a modified Horvitz-Thompson-type KDE to estimate the global density under missing data. However, several methods of local density estimation was widely discussed and extended in many directions; see~\citet{She83,She86},~\citet{Tho92},~\citet{Haz99},~\citet{Cha10} and recently~\citet{Dut14}.

In this paper, we propose a Recursive KDE to estimate the global unknown density $f$ of a univariate variable when the data are missing at random, we also give some comparative elements in the case of local density estimation. Then, we present a data-driven bandwidth selection algorithm for the proposed estimators. Data-driven bandwidth selection procedure was proposed by~\citet{Sla14a} in the framework of the recursive kernel density estimators, then,~\citet{Sla14b} propose a plug-in selection method for recursive kernel distribution estimators,~\citet{Sla15} propose a plug-in selection method for recursive kernel regression estimators with a fixed design setting, and~\citet{Sla16} propose a plug-in selection algorithm for the semi-recursive kernel regression estimators, all of this works suppose that the full data are observed. Here, we developed a specific second-generation plug-in bandwidth selection method of the Recursive KDE under missing data.\\

Interestingly, unlike the EM-algorithm, the proposed Recursive KDE has explicit form, which facilitates its usage and analysis. Moreover, the \texttt{CPU} time using the proposed recursive KDE are approximately two times faster than the \texttt{CPU} time using the nonrecursive Horvitz-Thompson-type KDE.

For further motivation in order to help intuition, let us introduce in some detail two datasets, the first one concerne a simulated data of the Aquitaine cohort of $\texttt{HIV}$-1 infected patients, see the paper of~\citet{Thi00}, the second one concerne a real data which correspond to two array Comparative Genomic Hybridization (CGH) studies of fibroblast cell strains. 
Assuming that all these measurement are corrects, we simulated various proportions of missing at random by removing the corresponding entries before performing our algorithm. The estimated densities are then compared to
 the shape of the histograms obtained using the full data. Even when $30\%$ of the original measurements are missing, the estimated densities using the proposed estimators remain very accurate thus demontrating the effectiveness of our approach. Moreover, we compared our approach with the nonrecursive Horvitz-Thompson-type KDE via real dataset as well as simulations.
Results showed that our approach outperformed other approaches in terms of estimation accuracy and computing efficiency.
The remainder of the paper is organized as follows. In the next section we present our proposed recursive KDE. In Section~\ref{section:2}, we state our main results. Section~\ref{section:app} is devoted to our application results, first by simulation (subsection~\ref{subsection:simu}) and second using two datasets, the first one concerne the simulated data of the Aquitaine cohort of $\texttt{HIV}$-1 infected patients and the second one concerne the coriell cell lines using the chromosome number 11 (subsection~\ref{subsection:real}). We conclude the paper in Section~\ref{section:conclusion}, wehereas the technical details are deferred to Appendix~\ref{section:proofs} and Appendix~\ref{section:results:HTKDE} gives the results for the Horvitz-Thompson-type KDE. 
\section{The estimator}\label{section:estimator}
Let $T,T_1,\ldots, T_n$ be independent, identically distributed random variables, and let $f$ and $F$ denote respectively the probability density and the distribution of $T$. We suppose that the full data $T_1,\ldots, T_n$ are not totally available and are subjects to missing data. The observed random variables are then $X_i$ and $\delta_i$ where
\begin{eqnarray*}
X_i=\delta_i\times T_i \,\mbox{and}\quad \delta_i=\mathds{1}_{\left\{T_i\,\mbox{is observed}\right\}},\,1\leq i\leq n.
\end{eqnarray*}
When some $T_i$ are missing, we let : $\pi_i=\mathbb{P}\left[\delta_i=1\vert T_i\right]\quad \mbox{for}\quad 1\leq i\leq n.$\\
This probability is often called the propensity score, see~\citet{Ros83}. Let us recall that, in order to construct a stochastic algorithm, which approximates the function $f$ at a given point $x$, we need to define an algorithm of search of the zero of the function $h : y\to f\left(x\right)-y$. Following Robbins-Monro's procedure, this algorithm is defined by setting $f_0(x)\in \mathbb{R}$, and, for all $n\geq 1$, 
\begin{eqnarray*}
f_n\left(x\right)=f_{n-1}\left(x\right)+\gamma_nW_n,
\end{eqnarray*}
where $W_n(x)$ is an observation of the function $h$ at the point $f_{n-1}(x)$, and the stepsize $\left(\gamma_n\right)$ is a sequence of positive real numbers that goes to zero. In order to define $W_n(x)$, we follow the approach of \citet{Rev73,Rev77}, \citet{Tsy90} and of \citet{Sla13,Sla14a,Sla14b}, and we introduce a kernel $K$ (that is, a function satisfying $\int_{\mathbb{R}} K(x)dx=1$), and a bandwidth $(h_n)$ (that is, a sequence of positive real numbers that goes to zero), and sets $W_n(x)=h_n^{-1}\delta_n\pi_n^{-1}K\left(h_n^{-1}\left[x-X_n\right]\right)-f_{n-1}(x)$.
Then, the estimator $f_n$ to recursively estimate the function $f$ at the point $x$ can be written as
\begin{eqnarray}\label{eq:fn}
f_n\left(x\right)=\left(1-\gamma_n\right)f_{n-1}\left(x\right)+\gamma_n\delta_n\pi_n^{-1}h_n^{-1}K\left(h_n^{-1}\left[x-X_n\right]\right),
\end{eqnarray}
where the stepsize $\left(\gamma_n\right)$ is a sequence of positive real numbers that goes to zero. Let us underline that we consider  $f_0\left(x\right)=0$ and we let $Q_n=\prod_{j=1}^n\left(1-\gamma_j\right)$, then the equation~(\ref{eq:fn}) can be rewritten as follows:
\begin{eqnarray*}
f_n\left(x\right)
&=&Q_n\sum_{k=1}^nQ_k^{-1}\gamma_k\delta_k\pi_k^{-1}h_k^{-1}K\left(\frac{x-X_k}{h_k}\right),
\end{eqnarray*}
which means that, we can estimate $f$ recursively at the point $x$ by using one of the two previous equations.\\  

Moreover, we show in the next section that the optimal bandwidth which minimize the $\mathbb{E}\int_{\mathbb{R}}\left[f_n\left(x\right)-f\left(x\right)\right]^2dx$  depend on the choice of the stepsizes $\left(\gamma_n\right)$; we show in particular that under some conditions of regularity of $f$ and using the stepsizes $\left(\gamma_n\right)=\left(\gamma_0\,n^{-1}\right)$, with $\gamma_0>2/5$, the bandwidth $\left(h_n\right)$ must equal 
\begin{eqnarray*}
\left(2^{-1/5}\left(\gamma_0-2/5\right)^{1/5}\left\{\frac{\int_{\mathbb{R}}f^2\left(x\right)dx}{\int_{\mathbb{R}}\left(f^{\left(2\right)}\left(x\right)f\left(x\right)\right)^2dx}\right\}^{1/5}\left\{\frac{\int_{\mathbb{R}}K^2\left(z\right)dz}{\left(\int_{\mathbb{R}}z^2K\left(z\right)dz\right)^2}\right\}^{1/5}\pi_n^{-1/5}n^{-1/5}\right).
\end{eqnarray*}
The first aim of the next section is to propose a data-driven bandwidth selection of the proposed recursive KDE in the case of missing data, and the second aim is to give the conditions under which the proposed estimators $f_n$ may behave more efficiently than the nonrecursive Horvitz-Thompson-type KDE see the papers~\citet{Hor52} and~\citet{Dub09}, defined as
\begin{eqnarray}\label{eq:Nad}
\widetilde{f}_n\left(x\right)=\frac{1}{nh_n}\sum_{i=1}^n\delta_i\pi_i^{-1}K\left(\frac{x-X_i}{h_n}\right).
\end{eqnarray}
In the next section we state our main results: the global density estimation is considered in Section~\ref{sub:sec:glo} and the local density estimation is developed in Section~\ref{sub:sec:loc}. 
 
\section{Assumptions and main results} \label{section:2}
We define the following class of regularly varying sequences.
\begin{defi}
Let $\gamma \in \mathbb{R} $ and $\left(v_n\right)_{n\geq 1}$ be a nonrandom positive sequence. We say that $\left(v_n\right) \in \mathcal{GS}\left(\gamma \right)$ if
\begin{eqnarray}\label{eq:5}
\lim_{n \to +\infty} n\left[1-\frac{v_{n-1}}{v_{n}}\right]=\gamma .
\end{eqnarray}
\end{defi}
Condition~(\ref{eq:5}) was introduced~\citet{Gal73} to define regularly varying sequences, see also the paper~\citet{Boj73} and by~\citet{Mok07} in the context of stochastic approximation algorithms. Noting that the acronym $\mathcal{GS}$ stand for (Galambos and Seneta). Typical sequences in $\mathcal{GS}\left(\gamma \right)$ are, for $b\in \mathbb{R}$, $n^{\gamma}\left(\log n\right)^{b}$, $n^{\gamma}\left(\log \log n\right)^{b}$, and so on. \\
In this section, we investigate asymptotic properties of the proposed estimators~(\ref{eq:fn}). The assumptions to which we shall refer are the following

\begin{description}
\item(A1) $K:\mathbb{R}\rightarrow \mathbb{R}$ is a continuous, bounded function satisfying $\int_{\mathbb{R}}K\left( z\right) dz=1$, and, $\int_{\mathbb{R}}zK\left(z\right)=0$ and $\int_{\mathbb{R}}z^2K\left(z\right)<\infty$. 
\item(A2) $i)$ $\left(\gamma_n\right)\in \mathcal{GS}\left(-\alpha\right)$ with $\alpha\in \left]1/2,1\right]$. \\
  $ii)$ $\left(h_n\right)\in \mathcal{GS}\left(-a\right)$ with $a\in \left]0,1\right[$.\\
  $iii)$ $\lim_{n\to \infty} \left(n\gamma_n\right)\in \left]\min\left\{2a,\left(\alpha-a\right)/2\right\},\infty\right]$.
\item(A3) $f$ is bounded, differentiable, and $f^{\left(2\right)}$ is bounded.
\end{description}
Assumption $\left(A2\right)(iii)$ on the limit of $\left(n\gamma_n\right)$ as $n$ goes to infinity is usual in the framework of stochastic approximation algorithms. It implies in particular that the limit of $\left(\left[n\gamma_n\right]^{-1}\right)$ is finite. For simplicity, we introduce the following notations:
\begin{eqnarray}
\xi&=&\lim_{n\to \infty}\left(n\gamma_n\right)^{-1},\label{eq:xi}\\
R\left(K\right)&=&\int_{\mathbb{R}}K^2\left(z\right)dz,\nonumber\\
\mu_j\left(K\right)&=&\int_{\mathbb{R}}z^jK\left(z\right)dz,\nonumber\\
\Theta\left(K\right)&=&R\left(K\right)^{4/5}\mu_2\left(K\right)^{2/5},\nonumber\\
I_1&=&\int_{\mathbb{R}}f^2\left(x\right)dx,\nonumber\\
I_2&=&\int_{\mathbb{R}}\left(f^{\left(2\right)}\left(x\right)\right)^2f\left(x\right)dx.\nonumber
\end{eqnarray} 
In this section, we explicit the choice of $\left(h_n\right)$ through a second-generation plug-in method, which minimize the Mean Weighted Integrated Squared Error $MWISE$ of the proposed recursive estimators~(\ref{eq:fn}) in the case of global estimation and through a second-generation plug-in method, which minimize the Mean Squared Error $MSE$ of the proposed recursive estimators~(\ref{eq:fn}) in the case of local estimation, in order to provide a comparison with the nonrecursive estimator~(\ref{eq:Nad}). Our first result is the following proposition, which gives the bias and the variance of $f_n$.
\begin{prop}[Bias and variance of $f_n$]\label{prop:bias:var:lambn}
Let Assumptions $\left(A1\right)-\left(A3\right)$ hold
\begin{enumerate}
\item If $a\in ]0, \alpha/5]$, then
\begin{eqnarray}\label{bias:lambn}
\mathbb{E}\left[f_n\left(x\right)\right]-f\left(x\right)=\frac{h_n^2}{2\left(1-2a\xi\right)}f^{\left(2\right)}\left(x\right)\mu_2\left(K\right)+o\left(h_n^2\right).
\end{eqnarray}
If $a\in  ]\alpha/5, 1[$, then
\begin{eqnarray}\label{bias:lambn:bis}
\mathbb{E}\left[f_n\left(x\right)\right]-f\left(x\right)=
o\left(\sqrt{\gamma_nh_n^{-1}}\right).
\end{eqnarray}
\item If $a\in [\alpha/5, 1[$, then
\begin{eqnarray}
Var\left[f_n\left(x\right)\right]=\frac{\gamma_n}{h_n}\pi_n^{-1}\frac{1}{\left(2-\left(\alpha-a\right)\xi\right)}f\left(x\right)R\left(K\right)+o\left(\frac{\gamma_n}{h_n}\right).
\label{var:lambn}
\end{eqnarray}
If $a\in ]0,\alpha/5[$, then
\begin{eqnarray}
Var\left[f_n\left(x\right)\right]=o\left(h_n^4\right).
\end{eqnarray}\label{var:lambn:bis}
\item If $\lim_{n\to \infty}\left(n\gamma_n\right)>\max\left\{2a, \left(\alpha-a\right)/2\right\}$, then~(\ref{bias:lambn}) and~(\ref{var:lambn}) hold simultaneously.
\end{enumerate}
\end{prop}
The bias and the variance of the proposed estimator $f_n$ defined by the stochastic approximation algorithm~(\ref{eq:fn}) then heavily depend on the choice of the stepsize $\left(\gamma_n\right)$. Let us first underline that, it follows from~(\ref{var:lambn}) that the stepsize which minimize the variance of $f_n$ is $\left(\gamma_n\right)=\left(\left[1-a\right]n^{-1}\right)$, using this stepsize the variance of $f_n$ is equal to
\begin{eqnarray*}
Var\left[f_n\left(x\right)\right]=\pi_n^{-1}\frac{1-a}{nh_n}f\left(x\right)R\left(K\right)+o\left(\frac{1}{nh_n}\right).
\end{eqnarray*}  
Now, using the special stepsize $\left(\gamma_n\right)=\left(n^{-1}\right)$,  (see~\citet{Mok09} and~\citet{Sla13}), the variance of $f_n$ is equal to
\begin{eqnarray*}
Var\left[f_n\left(x\right)\right]=\frac{\pi_n^{-1}}{1+a}\frac{1}{nh_n}f\left(x\right)R\left(K\right)+o\left(\frac{1}{nh_n}\right).
\end{eqnarray*} 
Let us recall that under the Assumptions $\left(A1\right)$, $\left(A2\right)ii)$ and $\left(A3\right)$, we have
\begin{eqnarray*}
Var\left[\widetilde{f}_n\left(x\right)\right]&=&\frac{\pi_n^{-1}}{nh_n}f\left(x\right)R\left(K\right)+o\left(\frac{1}{nh_n}\right).
\end{eqnarray*} 
Which shows that performing the proposed recursive estimators~(\ref{eq:fn}) with one of two proposed stepsizes we get smaller variance than the nonrecursive estimator~(\ref{eq:Nad}). Similar results was given by~\citet{Mok09} and~\citet{Sla13} in the case of complet data.

Let us now state the following theorem, which gives the weak convergence rate of the estimator $f_n$ defined in~(\ref{eq:fn}).
\begin{theor}[Weak pointwise convergence rate]\label{theo:TLC1}
Let Assumptions $\left(A1\right)-\left(A3\right)$ hold
\begin{enumerate}
\item If there exists $c\geq 0$ such that $\gamma_n^{-1}h_n^5\to c$, then
\begin{eqnarray*}
\lefteqn{\sqrt{\gamma_n^{-1}\pi_nh_n}\left(f_{n}\left(x\right)-f\left(x\right) \right) }\\
&&\stackrel{\mathcal{D}}{\rightarrow} 
\mathcal{N}\left(\frac{\sqrt{c}}{2\left(1-2a\xi\right)}f^{\left(2\right)}\left(x\right)\mu_2\left(K\right),\frac{1}{\left(2-\left(\alpha-a\right)\xi\right)}f\left(x\right)R\left(K\right)\right).
\end{eqnarray*}
\item If $\gamma_n^{-1}h_{n}^{5} \rightarrow \infty $, then  
\begin{eqnarray*}
\frac{1}{h_{n}^{2}}\left(f_{n}\left(x\right)-f\left(x\right)\right) \stackrel{\mathbb{P}}{\rightarrow } \frac{1}{2\left(1-2a\xi\right)}f^{\left(2\right)}\left(x\right)\mu_2\left(K\right),
\end{eqnarray*}
\end{enumerate}
where $\stackrel{\mathcal{D}}{\rightarrow}$ denotes the convergence in distribution, $\mathcal{N}$ the Gaussian-distribution and $\stackrel{\mathbb{P}}{\rightarrow}$ the convergence in probability.
\end{theor}
\subsection{Global density estimation}\label{sub:sec:glo}
In order to measure globally the quality of our recursive estimators~(\ref{eq:fn}), we use the following quantity, 
\begin{eqnarray*}
MWISE\left[f_n\right]&=&\mathbb{E}\int_{\mathbb{R}}\left[f_n\left(x\right)-f\left(x\right)\right]^2f\left(x\right)dx\nonumber\\
&=&\int_{\mathbb{R}}\left(\mathbb{E}\left(f_n\left(x\right)\right)-f\left(x\right)\right)^2f\left(x\right)dx+\int_{\mathbb{R}}Var\left(f_n\left(x\right)\right)f\left(x\right)dx.
\end{eqnarray*}
The following proposition gives the $MWISE$ of the recursive estimators defined in~(\ref{eq:fn}).
\begin{prop}[$MWISE$ of $f_n$]\label{prop:MISE:lamb}
Let Assumptions $\left(A1\right)-\left(A3\right)$ hold.
\begin{enumerate}
\item If $a\in  ]0, \alpha/5[$, then
\begin{eqnarray*}
MWISE\left[f_n\right]=\frac{h_n^4}{4}\frac{1}{\left(1-2a\xi\right)^2}I_2\mu_2^2\left(K\right)+o\left(h_n^4\right).
\end{eqnarray*}
\item If $a=\alpha/5$, then
\begin{eqnarray*}
MWISE\left[f_n\right]&=&\frac{\gamma_n\pi_n^{-1}}{h_n}\frac{1}{\left(2-\left(\alpha-a\right)\xi\right)}I_1R\left(K\right)+\frac{h_n^4}{4}\frac{1}{\left(1-2a\xi\right)^2}I_2\mu_2^2\left(K\right)+o\left(h_n^4\right).
\end{eqnarray*}
\item If $a\in ]\alpha/5, 1[$, then
\begin{eqnarray*}
MWISE\left[f_n\right]&=&\frac{\gamma_n\pi_n^{-1}}{h_n}\frac{1}{\left(2-\left(\alpha-a\right)\xi\right)}I_1R\left(K\right)
+o\left(\frac{\gamma_n\pi_n^{-1}}{h_n}\right).
\end{eqnarray*}
\end{enumerate}
\end{prop}
\subsubsection{Estimating propensity scores}
In order to estimate the propensity scores $\pi_i$, we need to exploit the information contained in the auxiliary variables $T_i$, which are related to $X_i$. First, we can estimate $\pi$ by the empirical proportion based on the observed data~(\citet{Qi05}):
\begin{eqnarray*}
\widehat{\pi}_i=\widehat{\pi}\left(X_i\right)=\frac{\sum_{j=1}^n\delta_j\mathds{1}_{\left\{X_j=X_i\right\}}}{\sum_{j=1}^n\mathds{1}_{\left\{T_j=T_i\right\}}}.
\end{eqnarray*}
Moreover, in the case when $X_i$ contains continuous elements,~\citet{Dub09} used the Nadaraya-Watson (local mean) estimator~(\citet{Nad64};~\citet{Wat64}):
\begin{eqnarray}\label{ps:NW}
\widehat{\pi}_{NWi}=\widehat{\pi}_{NW}\left(X_i\right)=\frac{\sum_{j=1}^n\delta_jK\left(h_n^{-1}\left[X_i-X_j\right]\right)}{\sum_{j=1}^nK\left(h_n^{-1}\left[X_i-X_j\right]\right)},
\end{eqnarray}
and in the case of binary nature of the response variable,~\citet{Tib87} and~\citet{Fan95} provide an estimation of the propensity score using local likelihood. However,~\citet{Dub09} concluded that, both the Nadaraya-Watson estimates and local likelihood estimates of the propensity scores are more flexible than ordinary binary regression estimates. In this paper we use the recursive version of the estimator of Nadaraya-Watson (called also the semi-recursive estimator) and defined as:
\begin{eqnarray}\label{ps:RNW}
\widehat{\pi}_{RNWi}=\widehat{\pi}_{RNW}\left(X_i\right)=\frac{\sum_{j=1}^n\delta_jh_j^{-1}K\left(h_j^{-1}\left[X_i-X_j\right]\right)}{\sum_{j=1}^nh_j^{-1}K\left(h_i^{-1}\left[X_i-X_j\right]\right)}.
\end{eqnarray}
For now, we simply use the propensity score estimator~(\ref{ps:NW}) in the case of nonrecursive Horvitz-Thompson-type KDE and~(\ref{ps:RNW}) in the case of the proposed recursive KDE.
\subsubsection{Bandwidth selection}
In the framework of the nonparametric kernel estimators, the bandwidth selection methods studied in the literature can be divided into three broad classes: the cross-validation techniques, the plug-in ideas and the bootstrap. A detailed comparison of the three practical bandwidth selection can be found in \citet{Del04}. They concluded that chosen appropriately plug-in and bootstrap selectors both outperform the cross-validation bandwidth, and that none of the two can be claimed to be best in all cases. In this section, we developed a plug-in bandwidth selector that minimizing the $MWISE$ of the proposed recursive KDE, using the function $f\left(x\right)$ as a weight function.\\
The following corollary ensures that the bandwidth which minimize the $MWISE$ of $f_n$ depend on the stepsize $\left(\gamma_n\right)$ and then the corresponding $MWISE$ depend also on the stepsize $\left(\gamma_n\right)$.
\begin{coro}\label{Coro:hn:MISE}
Let Assumptions $\left(A1\right)-\left(A3\right)$ hold. To minimize the $MWISE$ of $f_n$, the stepsize $\left(\gamma_n\right)$ must be chosen in $\mathcal{GS}\left(-1\right)$, the bandwidth $\left(h_n\right)$ must equal 
\begin{eqnarray*}
\left(\left\{\frac{\left(1-2a\xi\right)^2}{\left(2-\left(\alpha-a\right)\xi\right)}\frac{I_1}{I_2}\right\}^{1/5}\left\{\frac{R\left(K\right)}{\mu_2^2\left(K\right)}\right\}^{1/5}\pi_n^{-1/5}\gamma_n^{1/5}\right).
\end{eqnarray*}
Then, we have
\begin{eqnarray*} MWISE\left[f_n\right]&=&\frac{5}{4}\left(1-2a\xi\right)^{-2/5}\left(2-\left(\alpha-a\right)\xi\right)^{-4/5}I_1^{4/5}I_2^{1/5}\Theta\left(K\right)\pi_n^{-4/5}\gamma_n^{4/5}+o\left(\pi_n^{-4/5}\gamma_n^{4/5}\right).
\end{eqnarray*}
\end{coro}
The following corollary shows that, for a special choice of the stepsize $\left(\gamma_n\right)=\left(\gamma_0n^{-1}\right)$, which fulfilled that $\lim_{n\to \infty}n\gamma_n=\gamma_0$ and that $\left(\gamma_n\right)\in \mathcal{GS}\left(-1\right)$, the optimal value for $h_n$ depend on $\gamma_0$ and then the corresponding $MWISE$ depend on $\gamma_0$.
\begin{coro}\label{Coro:hn:MISE:bis}
Let Assumptions $\left(A1\right)-\left(A3\right)$ hold, and suppose that $\left(\gamma_n\right)=\left(\gamma_0n^{-1}\right)$. To minimize the $MWISE$ of $f_n$, the stepsize $\left(\gamma_n\right)$ must be chosen in $\mathcal{GS}\left(-1\right)$, $\lim_{n\to \infty}n\gamma_n=\gamma_0$, the bandwidth $\left(h_n\right)$ must equal 
\begin{eqnarray*}
\left(2^{-1/5}\left(\gamma_0-2/5\right)^{1/5}\left(\frac{I_1}{I_2}\right)^{1/5}\left\{\frac{R\left(K\right)}{\mu_2^2\left(K\right)}\right\}^{1/5}\pi_n^{-1/5}n^{-1/5}\right),
\end{eqnarray*}
and we then have
\begin{eqnarray*}
MWISE\left[f_n\right]&=&\frac{5}{4}\frac{1}{2^{4/5}}\frac{\gamma_0^2}{\left(\gamma_0-2/5\right)^{6/5}}I_1^{4/5}I_2^{1/5}\Theta\left(K\right)\pi_n^{-4/5}n^{-4/5}+o\left(\pi_n^{-4/5}n^{-4/5}\right).
\end{eqnarray*}
Moreover, the minimum of $\gamma_0^2\left(\gamma_0-2/5\right)^{-6/5}$ is reached at $\gamma_0=1$, then the bandwidth $\left(h_n\right)$ must equal 
\begin{eqnarray}\label{hoptim:gamma}
\left(\left(\frac{3}{10}\right)^{1/5}\left(\frac{I_1}{I_2}\right)^{1/5}\left\{\frac{R\left(K\right)}{\mu_2^2\left(K\right)}\right\}^{1/5}\pi_n^{-1/5}n^{-1/5}\right),
\end{eqnarray}
and then the asymptotic $MWISE$
\begin{eqnarray*} 
AMWISE\left[f_n\right]
&=&\frac{5}{4}\frac{1}{2^{4/5}}\left(\frac{5}{3}\right)^{6/5}I_1^{4/5}I_2^{1/5}\Theta\left(K\right)\pi_n^{-4/5}n^{-4/5}.
\end{eqnarray*}
\end{coro}
In order to estimate the optimal bandwidth (\ref{hoptim:gamma}), we must estimate the unknown quantites $I_1$ and $I_2$. For this purpose, we used the following modified version of the kernel estimators introduced in~\citet{Sla14a}:

\begin{eqnarray}
\widehat{I}_1&=&\frac{\Psi_n}{n}\sum_{i,k=1}^n\Psi_k^{-1}\beta_kb_k^{-1}\delta_k\widehat{\pi}_{RNWk}^{-1}K_b\left(\frac{X_i-X_k}{b_k}\right)\label{I1:hat:rec}\\
\widehat{I}_2&=&\frac{\Phi_n^2}{n}\sum_{\substack{i,j,k=1\\j\not=k}}^n\Phi_j^{-1}\Phi_k^{-1}\beta^{\prime}_j\beta^{\prime}_k\delta_j\delta_k\widehat{\pi}_{RNWj}^{-1}\widehat{\pi}_{RNWk}^{-1}b_j^{\prime-3}b_k^{\prime-3}K_{b^{\prime}}^{\left(2\right)}\left(\frac{X_i-X_j}{b_j^{\prime}}\right)K_{b^{\prime}}^{\left(2\right)}\left(\frac{X_i-X_k}{b_k^{\prime}}\right)\label{I2:hat:rec}
\end{eqnarray}
where $K_b$ and $K_{b^{\prime}}$ are a kernels, $b_n$ and $b^{\prime}_n$ are respectively the associated bandwidth (called pilot bandwidth) and $\beta_n$ and $\beta_n^{\prime}$ are the two pilot stepsizes for the estimation of $I_1$ and $I_2$ respectively, and $\Psi_n=\prod_{i=1}^n\left(1-\beta_i\right)$ and $\Phi_n=\prod_{i=1}^n\left(1-\beta_i^{\prime}\right)$.\\
In practice, we take
\begin{eqnarray}\label{eq:h:initial}
b_n=n^{-\beta}\min\left\{\widehat{s},\frac{Q_3-Q_1}{1.349}\right\},\quad\beta \in \left]0,1\right[
\end{eqnarray} 
(see~\citet{Sil86}) with $\widehat{s}$ the sample standard deviation, and $Q_1$, $Q_3$ denoting the first and third quartiles, respectively.\\
We followed the same steps as in~\citet{Sla14a} and we showed that in order to minimize the $AMISE$ of $\widehat{I}_1$ the pilot bandwidth $\left(b_n\right)$ must belong to $\mathcal{GS}\left(-2/5\right)$ and the pilot stepsize $\left(\beta_n\right)$ should be equal to $\left(1.36n^{-1}\right)$, and in order to minimize the $AMISE$ of $\widehat{I}_2$ the pilot bandwidth $\left(b_n^{\prime}\right)$ must belong to $\mathcal{GS}\left(-3/14\right)$ and the pilot stepsize $\left(\beta_n^{\prime}\right)$ should be equal to $\left(1.48n^{-1}\right)$.\\

Finally, the plug-in estimator of the bandwidth $\left(h_n\right)$ using the recursive estimators defined in~(\ref{eq:fn}) with the stepsizes $\left(\gamma_n\right)=\left(n^{-1}\right)$ is equal to 
\begin{eqnarray}\label{hoptim:gamma:MISE}
\left(\left(\frac{3}{10}\right)^{1/5}\left(\frac{\widehat{I}_1}{\widehat{I}_2}\right)^{1/5}\left\{\frac{R\left(K\right)}{\mu_2^2\left(K\right)}\right\}^{1/5}\widehat{\pi}_{RNWn}^{-1/5}n^{-1/5}\right),
\end{eqnarray}
and the associated plug-in $AMWISE$ is equal to
\begin{eqnarray*}
\widehat{AMWISE}\left[f_n\right]
&=&\frac{5}{4}\frac{1}{2^{4/5}}\left(\frac{5}{3}\right)^{6/5}\widehat{I}_1^{4/5}\widehat{I}_2^{1/5}\Theta\left(K\right)\widehat{\pi}_{RNWn}^{-4/5}n^{-4/5}.
\end{eqnarray*}
Now, using the stepsize $\left(\gamma_n\right)=\left(\left[4/5\right]n^{-1}\right)$, the stepsize which minimize the variance of the proposed estimators defined in~(\ref{eq:fn}), it follows from Corollary~\ref{Coro:hn:MISE:bis} that, the plug-in estimator of the bandwidth $\left(h_n\right)$ in this particular case is equal to
\begin{eqnarray}\label{hoptim:gamma:var}
\left(5^{-1/5}\left(\frac{\widehat{I}_1}{\widehat{I}_2}\right)^{1/5}\left\{\frac{R\left(K\right)}{\mu_2^2\left(K\right)}\right\}^{1/5}\widehat{\pi}_{RNWn}^{-1/5}n^{-1/5}\right),
\end{eqnarray}
and the associated plug-in $AMWISE$ is equal to
\begin{eqnarray*}
\widehat{AMWISE}\left[f_n\right]
&=&5^{1/5}\widehat{I}_1^{4/5}\widehat{I}_2^{1/5}\Theta\left(K\right)\widehat{\pi}_{RNWn}^{-4/5}n^{-4/5}.
\end{eqnarray*}
Moreover, following similar steps as in~\citet{Sla14a}, we prove the following corollary
\begin{coro}
Let the assumptions $\left(A1\right)-\left(A3\right)$ hold, and the bandwidth $\left(h_n\right)$ equal to~(\ref{hoptim:gamma:MISE}) and the stepsize $\left(\gamma_n\right)=\left(n^{-1}\right)$. We have
\begin{eqnarray*}
\frac{\mathbb{E}\left[\widehat{AMWISE}\left(f_n\right)\right]}{\mathbb{E}\left[\widetilde{AMWISE}\left(\widetilde{f}_n\right)\right]}<1. 
\end{eqnarray*}
Then, the expectation of the estimated $AMWISE$ of the proposed recursive estimators defined by~(\ref{eq:fn}) using the special stepsize $\left(\gamma_n\right)=\left(n^{-1}\right)$ is smaller than the expectation of the estimated $AMWISE$ of the nonrecursive Horvitz-Thompson-type KDE defined by~(\ref{eq:Nad}) (see Appendix~\ref{section:results:HTKDE} for some results on the nonrecursive Horvitz-Thompson-type KDE estimator.
\end{coro} 
\subsection{Local density estimation}\label{sub:sec:loc}
In order to measure locally the quality of our recursive estimators~(\ref{eq:fn}), we use the following quantity, 
\begin{eqnarray*}
MSE\left[f_n\right]&=&\mathbb{E}\left[f_n\left(x\right)-f\left(x\right)\right]^2\nonumber\\
&=&\left(\mathbb{E}\left(f_n\left(x\right)\right)-f\left(x\right)\right)^2+Var\left(f_n\left(x\right)\right).
\end{eqnarray*}
Following similar step as in the case of the global density estimation, we show that the value of $h_n$ which minimizes the asymptotic mean square error of $f_n\left(x\right)$ is equal to
\begin{eqnarray}\label{hoptim:gamma:loc}
\left(\left(\frac{3}{10}\right)^{1/5}\left(\frac{f\left(x\right)}{\left(f^{\left(2\right)}\left(x\right)\right)^2}\right)^{1/5}\left\{\frac{R\left(K\right)}{\mu_2^2\left(K\right)}\right\}^{1/5}\pi_{n}^{-1/5}n^{-1/5}\right).
\end{eqnarray}
In practice the kernel $K$ is known but $f\left(x\right)$ and $f^{\left(2\right)}\left(x\right)$ are not. Thus, we estimate $f\left(x\right)$ by $f_n\left(x\right)$ and $f^{\left(2\right)}\left(x\right)$ by $f_n^{\left(2\right)}\left(x\right)$, and so the plug-in bandwidth selection in the case of locally estimation is equal to
\begin{eqnarray}\label{hoptim:gamma:loc:plug:in}
\left(\left(\frac{3}{10}\right)^{1/5}\left(\frac{\widehat{f}_n\left(x\right)}{\left(\widehat{f}_n^{\left(2\right)}\left(x\right)\right)^2}\right)^{1/5}\left\{\frac{R\left(K\right)}{\mu_2^2\left(K\right)}\right\}^{1/5}\widehat{\pi}_{RNWn}^{-1/5}n^{-1/5}\right),
\end{eqnarray}
where 
\begin{eqnarray*}
\widehat{f}_n\left(x\right)&=&\Psi_n\sum_{k=1}^n\Psi_k^{-1}\beta_kb_k^{-1}\delta_k\widehat{\pi}_{RNWk}^{-1}K_b\left(\frac{x-X_k}{b_k}\right)\label{fn:hat:rec:loc}\\
\widehat{f}^{\left(2\right)}_n\left(x\right)&=&\Phi_n\sum_{k=1}^n\Phi_k^{-1}\beta^{\prime}_k\delta_k\widehat{\pi}_{RNWk}^{-1}b_k^{\prime-3}K_{b^{\prime}}^{\left(2\right)}\left(\frac{x-X_k}{b_k^{\prime}}\right)\label{fn:sec:hat:rec:loc}
\end{eqnarray*}
where $K_b$ and $K_{b^{\prime}}$ are a kernels, $b_n$ and $b^{\prime}_n$ are respectively the associated bandwidth (called pilot bandwidth) and $\beta_n$ and $\beta_n^{\prime}$ are the two pilot stepsizes for the estimation of $f\left(x\right)$ and $f^{\left(2\right)}\left(x\right)$ respectively, and $\Psi_n=\prod_{i=1}^n\left(1-\beta_i\right)$ and $\Phi_n=\prod_{i=1}^n\left(1-\beta_i^{\prime}\right)$.\\
\section{Applications}\label{section:app}
The aim of our applications is to compare the performance of the recursive kernel density estimators under missing data defined in~(\ref{eq:fn}) with that of the nonrecursive Horvitz-Thompson-type KDE defined in~(\ref{eq:Nad}). 
\begin{description}
\item When applying $f_n$ one need to choose four quantities:
\begin{enumerate}
\item The function $K$, we choose the Normal kernel. 
\item The stepsizes $\left(\gamma_n\right)$ equal to $\left(c\,n^{-1}\right)$, with $c\in \left[4/5,1\right]$.
\item The propensity score $\left(\pi_n\right)$ is chosen to be equal to~(\ref{ps:RNW}).
\item The bandwidth $\left(h_n\right)$ is chosen to be equal  to~(\ref{hoptim:gamma:MISE}) in the case when $\left(\gamma_n\right)=\left(n^{-1}\right)$ and to~(\ref{hoptim:gamma:var}) in the case when $\left(\gamma_n\right)=\left(\left[4/5\right]n^{-1}\right)$.
\begin{enumerate}
\item To estimate $I_1$, we use (\ref{I1:hat:rec}); The pilot bandwidth is chosen to be equal to~(\ref{eq:h:initial}) with the choice of $\beta=2/5$ and the pilot stepsize equal to $\left(1.36n^{-1}\right)$.
\item To estimate $I_2$, we use (\ref{I2:hat:rec}); The pilot bandwidth is chosen to be equal to~(\ref{eq:h:initial}) with the choice of $\beta=3/14$ and the pilot stepsize equal to $\left(1.48n^{-1}\right)$.
\end{enumerate}
\end{enumerate}
\item When applying $\widetilde{f}_n$ one need to choose three quantities: 
\begin{enumerate}
\item The function $K$, as in the recursive framework, we use the Normal kernel.
\item The propensity score $\left(\pi_n\right)$ is chosen to be equal to~(\ref{ps:NW}).
\item The bandwidth $\left(h_n\right)$ is chosen to be equal to~(\ref{hoptim:rose:plug:in}).
\begin{enumerate}
\item To estimate $I_1$, we use (\ref{I1:norec}); The pilot bandwidth is chosen to be equal to~(\ref{eq:h:initial}) with the choice of $\beta=2/5$.
\item To estimate $I_2$, we use (\ref{I2:norec}); The pilot bandwidth is chosen to be equal to~(\ref{eq:h:initial}) with the choice of $\beta=3/14$. 
\end{enumerate}
\end{enumerate}
\end{description}

\subsection{Simulations}\label{subsection:simu}
\subsubsection{Global density estimation}
In order to investigate the comparison between the proposed estimators, we consider three sample sizes: $n=100$, $200$, and $500$. In each case, we consider $0\%$, $30\%$, $50\%$ and $70\%$ of missing data, the number of the design points was fixed to be equal to $500$, which variate from the lowest value to highest value over ($N=500$ number of simulations) with equally spaced setting. We then use the bandwidth selection method proposed in the applications Section. Moreover, we consider three densities functions $f$ : 1- the standard normal : $X\sim \mathcal{N}\left(0,1\right)$ (see Table~\ref{Tab:1}), 2- the normal mixture distribution : $X\sim 1/2\mathcal{N}\left(2,1\right)+1/2\mathcal{N}\left(-3,1\right)$ (see Table~\ref{Tab:2}), 3- the weibull distribution with shape parameter $2$ and scale parameter $1$: $X\sim \mathcal{W}eibul\left(2,1\right)$ (see Table~\ref{Tab:3}). We considered also the case of one auxiliary variable and generated the pairs $\left(X,Y\right),\left(X_1,Y_1\right),\ldots, \left(X_n,Y_n\right)$ by first generating $Y$ distributed as $F_Y$ and $X^*$ distributed as $F_{X^*}$, and we let $X=\rho Y+\sqrt{1-\rho^2}X^*$, $-1<\rho<1$, then, $Y$ and $X$ were correlated. Moreover, we generated $X^*,X_1^*,\ldots,X_n^*$ from the standard normal distribution and we set $\rho=0.3$, $0.5$ and $0.8$ and we generated the responses $Y,Y_1,\ldots,Y_n$ from one of the three previous densities functions $f$ (see Table~\ref{Tab:4}), in this case we present results for the situations in which $70\%$ of the responses are missing. For each situations, we compute the Mean Weighted Integrated Squared Error $MWISE$ (over $500$ samples).

\paragraph{Computational cost}
The advantage of recursive estimators on their nonrecursive version is that their update, from a sample of size $n$ to one of size $n+1$, require less computations. This property can be generalized, one can check that it follows from~(\ref{eq:fn}) that for all $n_1\in \left[0,n-1\right]$,
\begin{eqnarray*}
f_n\left(x\right)=\prod_{j=n_1+1}^n\left(1-\gamma_j\right)f_{n_1}\left(x\right)+\sum_{k=n_1}^{n-1}\left[\prod_{j=k+1}^n\left(1-\gamma_j\right)\right]\gamma_kh_k^{-1}\delta_k\pi_k^{-1}K\left(h_k^{-1}\left(x-X_k\right)\right).
\end{eqnarray*}
In order to give some comparative elements with nonrecursive Horvitz-Thompson-type KDE~(\ref{eq:Nad}), including computational costs. We consider a $500$ samples of size $n_1=\lfloor n/2\rfloor$ (the lower integer part of $n/2$) generated from respectively the five considered distributions, moreover, we suppose that we receive an additional $500$ samples of size $n-n_1$ generated also from the same five considered distributions. 
Performing the two methods, we report the total \texttt{CPU} time values for each considered distribution and in all cases given in Tables~\ref{Tab:1},~\ref{Tab:2},~\ref{Tab:3} and~\ref{Tab:4}, the \texttt{CPU} time is given in seconds.

\begin{table}
\begin{center}
\begin{tabular}{lcc|cc|cc}
& \multicolumn{6}{c}{$X\sim\mathcal{N}\left(0,1\right)$}\\
& \multicolumn{2}{c|}{$n=100$}&\multicolumn{2}{c|}{$n=200$}&\multicolumn{2}{c}{$n=500$} \\
\hline
$0\%$&MWISE &CPU& MWISE &CPU& MWISE &CPU\\
\hline
Nonrecursive&$3.86e^{-03}$&$364$& $2.33e^{-03}$&$2448$& $1.20e^{-03}$& $14432$\\
Recursive 1& ${\bf 3.70e^{-03}}$&$194$& ${\bf 2.30e^{-03}}$&$1225$& $ {\bf 1.20e^{-03}}$&${\bf 7118}$\\
Recursive 2& $3.85e^{-03}$&${\bf 184}$& $2.39e^{-03}$&${\bf 1169}$& $ 1.25e^{-03}$&$7122$\\
\hline
$30\%$&MWISE &CPU& MWISE &CPU& MWISE &CPU\\
\hline
Nonrecursive&$1.47e^{-02}$&$365$& $5.93e^{-03}$&$2448$& $ 2.59e^{-03}$&$14134$\\
Recursive 1& ${\bf 1.27e^{-02}}$&${\bf 164}$& $ {\bf 5.73e^{-03}}$&${\bf 1206}$&$ {\bf 2.56e^{-03}}$&$7112$\\
Recursive 2& $ 1.32e^{-02}$&$167$& $ 6.17e^{-03}$&$1216$&$ 2.70e^{-03}$&${\bf 7043}$\\
\hline
$50\%$&MWISE &CPU& MWISE &CPU& MWISE &CPU\\
\hline
Nonrecursive&$2.23e^{-02}$&$356$& $1.87e^{-02}$&$2567$& $ 3.32e^{-03}$&$14123$\\
Recursive 1& ${\bf 2.11e^{-02}}$&${\bf 169}$& ${\bf 1.69e^{-02}}$&${\bf 1146}$& ${\bf 3.24e^{-03}}$&${\bf 7089}$\\
Recursive 2& $2.17e^{-02}$&$176$& $1.85e^{-02}$&$1154$& $3.46e^{-03}$&$7134$\\
\hline
$70\%$&MWISE &CPU& MWISE &CPU& MWISE &CPU\\
\hline
Nonrecursive&$2.89e^{-02}$&$346$& $2.43e^{-02}$&$2267$& $ 4.82e^{-03}$&$14134$\\
Recursive 1& ${\bf 2.63e^{-02}}$ &${\bf 168}$ &${\bf 2.32e^{-02}}$&${\bf 1166}$& ${\bf 4.63e^{-03}}$&$7156$\\
Recursive 2& $2.77e^{-02}$ &$177$ &$2.39e^{-02}$&$1175$& $4.86e^{-03}$&${\bf 7143}$\\
\hline
\end{tabular}
\end{center}
\caption{Quantitative comparison between the nonrecursive Horvitz-Thompson-type KDE~(\ref{eq:Nad}) and two recursive estimators; recursive $1$ correspond to the proposed estimator~(\ref{eq:fn}) with the choice of $\left(\gamma_n\right)=\left(n^{-1}\right)$, resursive $2$ correspond to the estimator~(\ref{eq:fn}) with the choice of $\left(\gamma_n\right)=\left(\left[4/5\right]n^{-1}\right)$. Here we consider the normal distribution $X\sim \mathcal{N}\left(0,1\right)$ with the proportion of missing at random equal respectively to $0\%$ (in the first block), equal to $10\%$ (in the second block), equal to $20\%$ (in the third block) and equal to $30\%$ (in the last block), we consider three sample sizes $n=100$, $n=200$ and $n=500$, the number of simulations is $500$, and we compute the Mean Weighted Integrated Squared Error ($MWISE$) and the \texttt{CPU} time in seconds.}\label{Tab:1}
\end{table}

\begin{table}
\begin{center}
\begin{tabular}{lcc|cc|cc}
& \multicolumn{6}{c}{$X\sim\frac{1}{2}\mathcal{N}\left(2,1\right)+\frac{1}{2}\mathcal{N}\left(-3,1\right))$}\\
& \multicolumn{2}{c|}{$n=100$}&\multicolumn{2}{c|}{$n=200$}&\multicolumn{2}{c}{$n=500$} \\
\hline
$0\%$&MWISE &CPU& MWISE &CPU& MWISE &CPU\\
\hline
Nonrecursive&$9.96e^{-04}$&$372$& $6.09e^{-04}$& $2448$& $3.38e^{-04}$&$14104$\\
Recursive 1& $8.64e^{-04}$&${\bf 174}$& ${\bf 5.50e^{-04}}$&${\bf 1168}$& $ {\bf 3.09e^{-04}}$&$7016$\\
Recursive 2& $ 9.00e^{-04}$&$177$& $5.72e^{-04}$&$1189$& $ 3.23e^{-04}$&${\bf 7003}$\\
\hline
$30\%$&MWISE &CPU& MWISE &CPU& MWISE &CPU\\
\hline
Nonrecursive&$1.30e^{-03}$&$366$& $8.31e^{-04}$&$2449$& $ 3.90e^{-04}$&$14142$\\
Recursive 1& ${\bf 1.25e^{-03}}$&$175$& $ {\bf 7.69e^{-04}}$&$1246$&$ {\bf 3.82e^{-04}}$&${\bf 7215}$\\
Recursive 2& $1.30e^{-03}$&${\bf 168}$& $ 8.01e^{-04}$&${\bf 1238}$&$ 3.97e^{-04}$&$7225$\\
\hline
$50\%$&MWISE &CPU& MWISE &CPU& MWISE &CPU\\
\hline
Nonrecursive&$1.41e^{-03}$&$354$& $1.55e^{-03}$&$2246$& $ 8.07e^{-04}$&$13994$\\
Recursive 1& ${\bf 1.29e^{-03}}$&${\bf 166}$& ${\bf 1.31e^{-03}}$&$1185$& ${\bf 6.52e^{-04}}$&${\bf 6946}$\\
Recursive 2& $1.35e^{-03}$&$173$& $1.36e^{-03}$&${\bf 1174}$& $6.79e^{-04}$&$6987$\\
\hline
$70\%$&MWISE &CPU& MWISE &CPU& MWISE &CPU\\
\hline
Nonrecursive&$1.95e^{-03}$&$356$& $1.87e^{-03}$&$2347$& $ 1.13e^{-03}$&$13965$\\
Recursive 1& ${\bf 1.80e^{-03}}$ &${\bf 168}$ &${\bf 1.69e^{-03}}$&$1164$& ${\bf 1.09e^{-03}}$&${\bf 7064}$\\
Recursive 2& $2.02e^{-03}$ &$179$ &$1.72e^{-03}$&${\bf 1162}$& $1.11e^{-03}$&$7086$\\
\hline
\end{tabular}
\end{center}
\caption{
Quantitative comparison between the nonrecursive Horvitz-Thompson-type KDE~(\ref{eq:Nad}) and two recursive estimators; recursive $1$ correspond to the proposed estimator~(\ref{eq:fn}) with the choice of $\left(\gamma_n\right)=\left(n^{-1}\right)$, resursive $2$ correspond to the estimator~(\ref{eq:fn}) with the choice of $\left(\gamma_n\right)=\left(\left[4/5\right]n^{-1}\right)$. Here we consider the normal mixture distribution $X\sim \frac{1}{2}\mathcal{N}\left(2,1\right)+\frac{1}{2}\mathcal{N}\left(-3,1\right)$, with the proportion of missing at random equal respectively to $0\%$ (in the first block), equal to $30\%$ (in the second block), equal to $50\%$ (in the third block) and equal to $70\%$ (in the last block), we consider three sample sizes $n=100$, $n=200$ and $n=500$, the number of simulations is $500$, and we compute the Mean Weighted Integrated Squared Error ($MWISE$) and the \texttt{CPU} time in seconds.}\label{Tab:2}
\end{table}

\begin{table}
\begin{center}
\begin{tabular}{lcc|cc|cc}
& \multicolumn{6}{c}{$X\sim \mathcal{W}eibull\left(2,1\right)$}\\
& \multicolumn{2}{c|}{$n=100$}&\multicolumn{2}{c|}{$n=200$}&\multicolumn{2}{c}{$n=500$} \\
\hline
$0\%$&MWISE &CPU& MWISE &CPU& MWISE &CPU\\
\hline
Nonrecursive&$1.82e^{-02}$&$374$& $1.11e^{-02}$& $2357$& $4.75e^{-03}$&$14256$\\
Recursive 1& $1.70e^{-02}$&${\bf 164}$& ${\bf 1.09e^{-02}}$&${\bf 1184}$& $ {\bf 4.74e^{-04}}$&$7122$\\
Recursive 2& $1.77e^{-04}$&$178$& $1.15e^{-02}$&$1194$& $ 4.94e^{-04}$&${\bf 7089}$\\
\hline
$30\%$&MWISE &CPU& MWISE &CPU& MWISE &CPU\\
\hline
Nonrecursive&$2.09e^{-02}$&$366$& $1.22e^{-02}$&$2452$& $ 5.19e^{-03}$&$14327$\\
Recursive 1& ${\bf 1.93e^{-03}}$&$177$& $ {\bf 1.16e^{-04}}$&$1154$&$ {\bf 5.12e^{-03}}$&$7012$\\
Recursive 2& $2.01e^{-03}$&${\bf 167}$& $ 1.23e^{-04}$&${\bf 1124}$&$ 5.27e^{-03}$&${\bf 6984}$\\
\hline
$50\%$&MWISE &CPU& MWISE &CPU& MWISE &CPU\\
\hline
Nonrecursive&$2.25e^{-02}$&$356$& $1.65e^{-02}$&$2354$& $ 1.14e^{-02}$&$14232$\\
Recursive 1& ${\bf 2.14e^{-03}}$&${\bf 166}$& ${\bf 1.59e^{-02}}$&${\bf 1163}$& ${\bf 1.07e^{-02}}$&$7048$\\
Recursive 2& $2.19e^{-02}$&$168$& $1.62e^{-02}$&$1165$& $1.12e^{-02}$&${\bf 7044}$\\
\hline
$70\%$&MWISE &CPU& MWISE &CPU& MWISE &CPU\\
\hline
Nonrecursive&$2.45e^{-02}$&$364$& $2.12e^{-02}$&$2346$& $ 1.82e^{-02}$&$13982$\\
Recursive 1& ${\bf 2.36e^{-02}}$ &$174$ &${\bf 2.03e^{-02}}$&$1187$& ${\bf 1.72e^{-02}}$&${\bf 7086}$\\
Recursive 2& $2.41e^{-02}$ &${\bf 166}$ &$2.08e^{-02}$&${\bf 1184}$& $1.83e^{-02}$ &$7094$\\
\hline
\end{tabular}
\end{center}
\caption{Quantitative comparison between the nonrecursive Horvitz-Thompson-type KDE~(\ref{eq:Nad}) and two recursive estimators; recursive $1$ correspond to the proposed estimator~(\ref{eq:fn}) with the choice of $\left(\gamma_n\right)=\left(n^{-1}\right)$, resursive $2$ correspond to the estimator~(\ref{eq:fn}) with the choice of $\left(\gamma_n\right)=\left(\left[4/5\right]n^{-1}\right)$. Here we consider the weibull distribution with shape parameter $2$ and scale parameter $1$, $X\sim \mathcal{W}eibull\left(2,1\right)$, with the proportion of missing at random equal respectively to $0\%$ (in the first block), equal to $30\%$ (in the second block), equal to $50\%$ (in the third block) and equal to $70\%$ (in the last block), we consider three sample sizes $n=100$, $n=200$ and $n=500$, the number of simulations is $500$, and we compute the Mean Weighted Integrated Squared Error $MWISE$ and the \texttt{CPU} time in seconds.}\label{Tab:3}
\end{table}

\begin{center}
\begin{figure}
    \includegraphics[width=0.6\textwidth,angle=270,clip=true,trim=40 0 0 0]{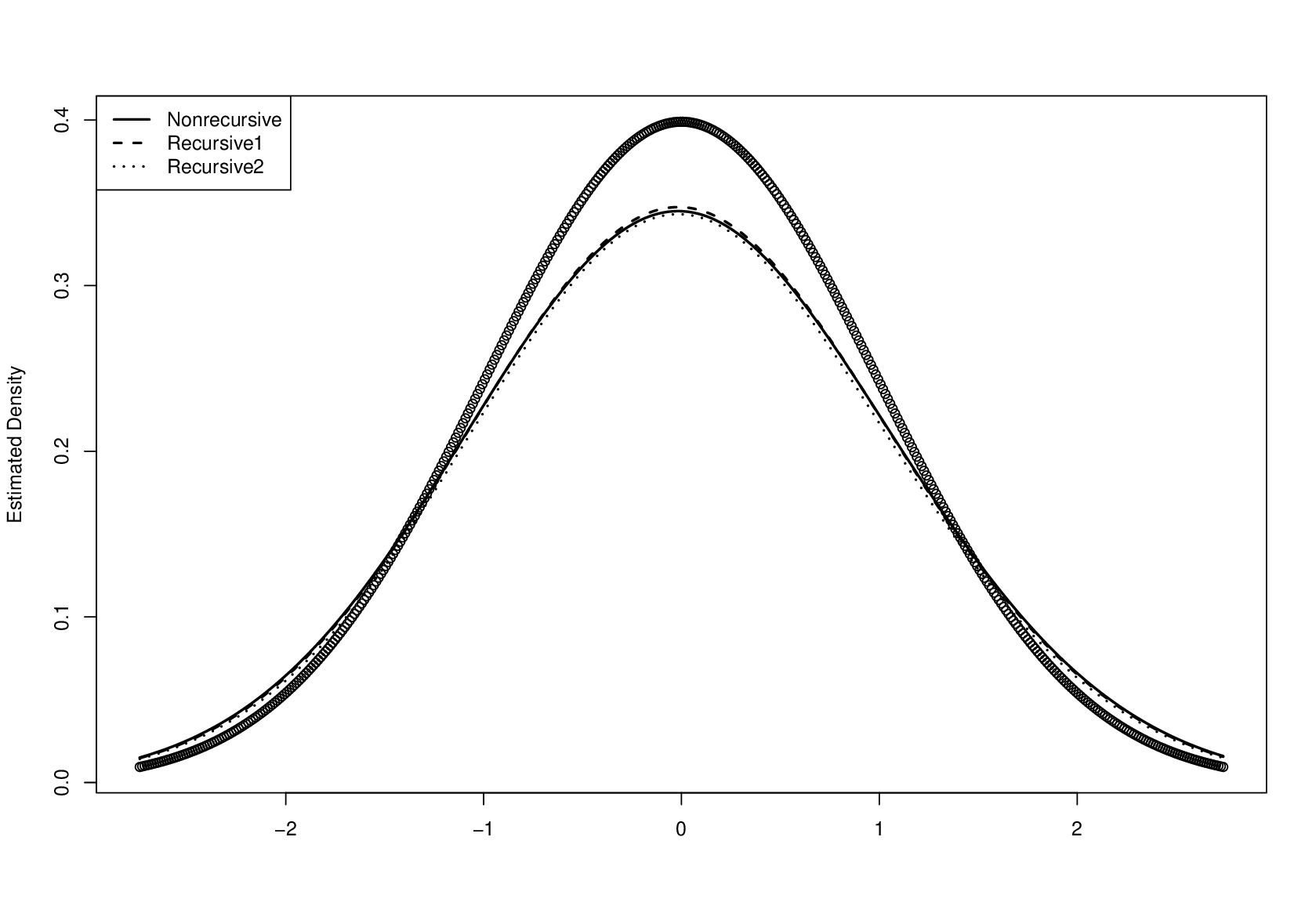}
\caption{Qualitative comparison between the nonrecursive Horvitz-Thompson-type KDE~(\ref{eq:Nad}) and two recursive estimators; recursive $1$ correspond to the proposed estimator~(\ref{eq:fn}) with the choice of $\left(\gamma_n\right)=\left(n^{-1}\right)$, resursive $2$ correspond to the estimator~(\ref{eq:fn}) with the choice of $\left(\gamma_n\right)=\left(\left[4/5\right]n^{-1}\right)$. Here we consider the normal distribution $\mathcal{N}\left(0,1\right)$ under $70\%$ of the missing data and we consider $500$ samples of size $n=200$.}
\label{Fig:1}
\end{figure}
\end{center}

\begin{center}
\begin{figure}
   \includegraphics[width=0.6\textwidth,angle=270,clip=true,trim=40 0 0 0]{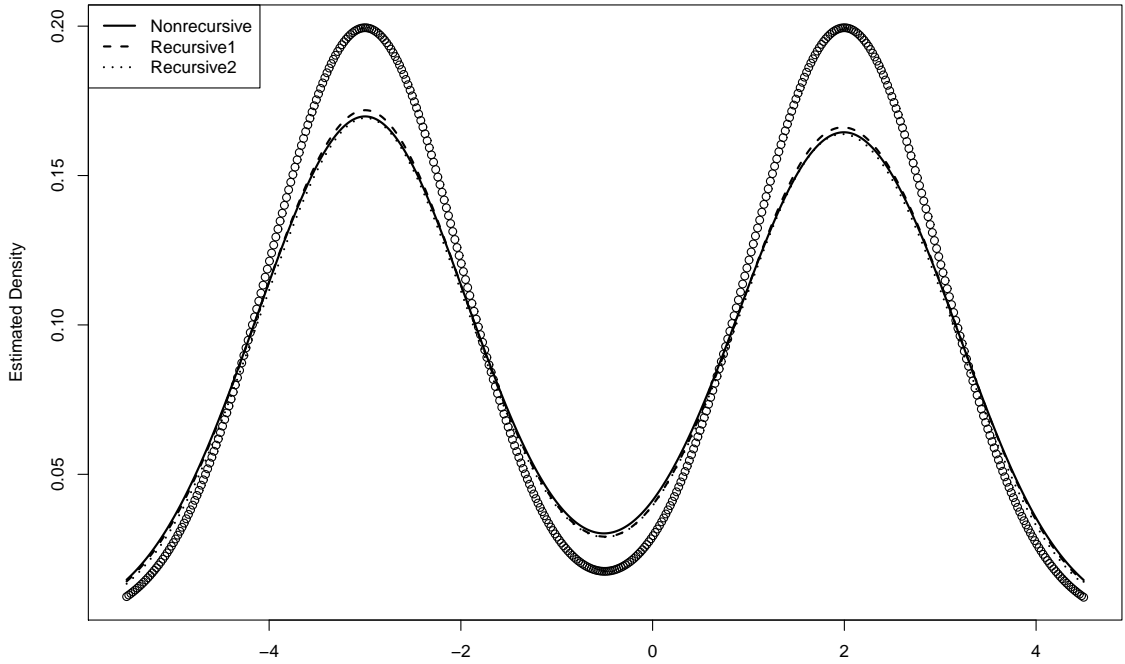}
\caption{Qualitative comparison between the nonrecursive Horvitz-Thompson-type KDE~(\ref{eq:Nad}) and two recursive estimators; recursive $1$ correspond to the proposed estimator~(\ref{eq:fn}) with the choice of $\left(\gamma_n\right)=\left(n^{-1}\right)$, resursive $2$ correspond to the estimator~(\ref{eq:fn}) with the choice of $\left(\gamma_n\right)=\left(\left[4/5\right]n^{-1}\right)$. Here we consider the normal mixture distribution $X\sim \frac{1}{2}\mathcal{N}\left(2,1\right)+\frac{1}{2}\mathcal{N}\left(-3,1\right)$  under $70\%$ of the missing data and we consider $500$ samples of size $n=200$.}
\label{Fig:2}
\end{figure}
\end{center}

\begin{center}
\begin{figure}
   \includegraphics[width=0.6\textwidth,angle=270,clip=true,trim=40 0 0 0]{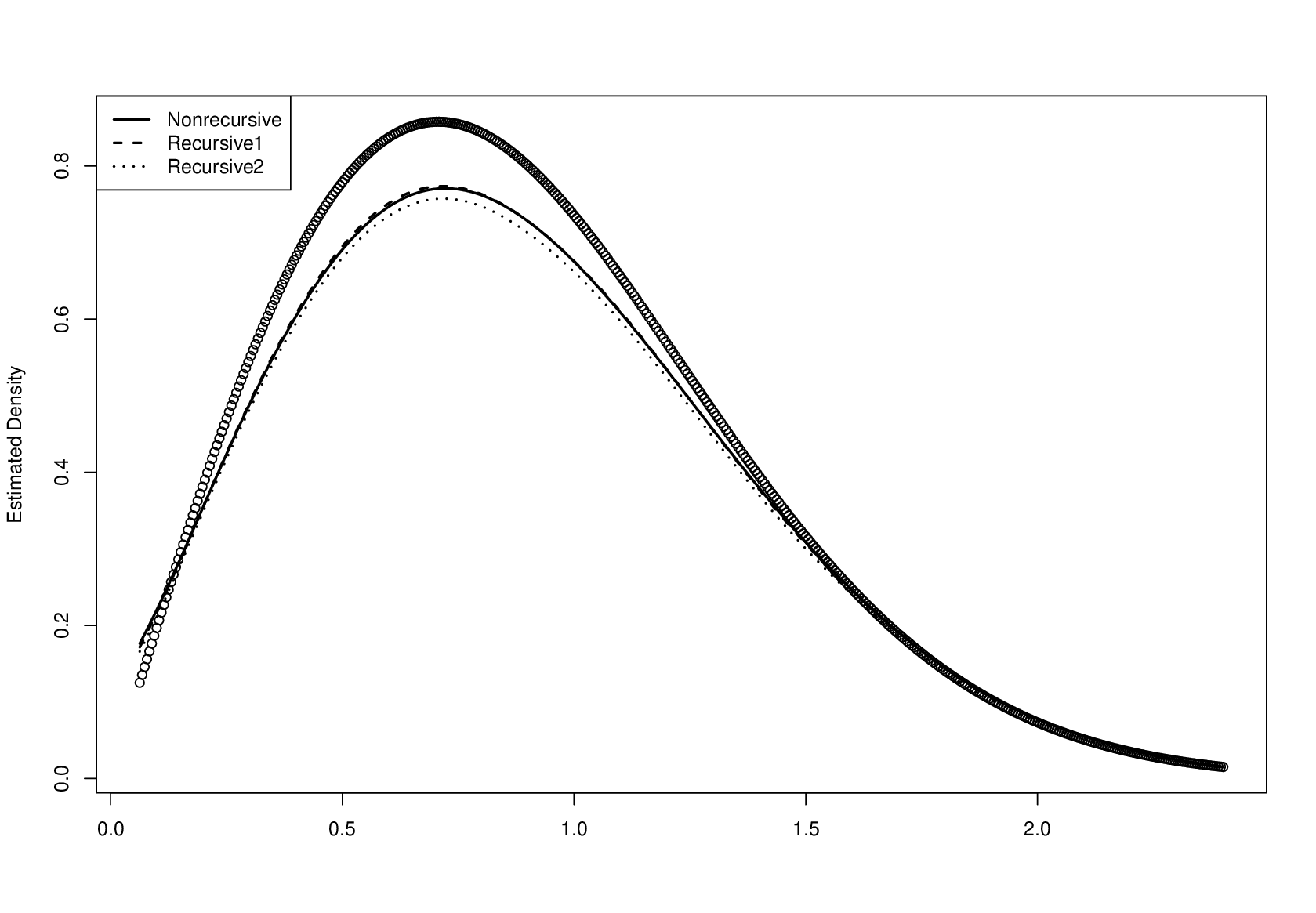}
\caption{Qualitative comparison between the nonrecursive Horvitz-Thompson-type KDE estimator~(\ref{eq:Nad}) and two recursive estimators; recursive $1$ correspond to the proposed estimator~(\ref{eq:fn}) with the choice of $\left(\gamma_n\right)=\left(n^{-1}\right)$, resursive $2$ correspond to the estimator~(\ref{eq:fn}) with the choice of $\left(\gamma_n\right)=\left(\left[4/5\right]n^{-1}\right)$. Here we consider the weibull distribution with shape parameter $2$ and scale parameter $1$, $X\sim \mathcal{W}eibull\left(2,1\right)$ under $70\%$ of the missing data and we consider $500$ samples of size $n=200$.}
\label{Fig:3}
\end{figure}
\end{center}

\begin{table}
\begin{center}
\begin{tabular}{lcc|cc|cc}
\hline
& \multicolumn{2}{c|}{$n=100$}&\multicolumn{2}{c|}{$n=200$}&\multicolumn{2}{c}{$n=500$} \\
\hline
& \multicolumn{6}{c}{$Y\sim \mathcal{N}\left(0,1\right)$}\\
\hline
$\rho=30\%$&MWISE &CPU& MWISE &CPU& MWISE &CPU\\
\hline
Nonrecursive&$0.02882$& $357$& $0.02584$& $2335$& $0.05042$&$14342$\\
Recursive 1& ${\bf 0.02622}$&${\bf 187}$& ${\bf 0.0244}$&${\bf 1244}$& $ {\bf 0.0484}$&${\bf 7332}$\\
Recursive 2& $0.02746$&$194$& $0.02542$&$1255$& $0.05044$&$7412$\\
\hline
$\rho=50\%$&MWISE &CPU& MWISE &CPU& MWISE &CPU\\
\hline
Nonrecursive&$0.02866$&$366$& $0.02564$&$2534$& $0.05054$&$14542$\\
Recursive 1& ${\bf 0.02634}$&$191$& $ {\bf 0.02454}$&$1254$&$ {\bf 0.04864}$&${\bf 7353}$\\
Recursive 2& $0.02724$&${\bf 187}$& $0.02544$&${\bf 1243}$&$ 0.05124$&$7366$\\
\hline
$\rho=80\%$&MWISE &CPU& MWISE &CPU& MWISE &CPU\\
\hline
Nonrecursive&$0.02894$&$362$& $0.02586$&$2447$& $0.05084$&$14121$\\
Recursive 1& ${\bf 0.02664}$& $188$& $ {\bf 0.02486}$&$1214$&$ {\bf 0.00493}$& ${\bf 7215}$\\
Recursive 2& $0.02784$& ${\bf 184}$& $0.02564$&${\bf 1209}$&$ 0.05164$&$7226$\\
\hline
& \multicolumn{6}{c}{$Y\sim 1/2\mathcal{N}\left(2,1\right)+1/2\mathcal{N}\left(-3,1\right)$}\\
\hline
$\rho=30\%$&MWISE &CPU& MWISE &CPU& MWISE &CPU\\
\hline
Nonrecursive&$0.00328$& $354$& $0.00314$&$2126$& $0.00186$&$14054$\\
Recursive 1& ${\bf 0.03014}$&$184$& $ {\bf 0.00306}$& $1207$&$ {\bf 0.00172}$&${\bf 7105}$\\
Recursive 2& $0.00329$& $ {\bf 177}$& $0.00316$& ${\bf 1198}$&$ 0.00187$&$7115$\\
\hline
$\rho=50\%$&MWISE &CPU& MWISE &CPU& MWISE &CPU\\
\hline
Nonrecursive&$0.00332$&$357$& $0.00221$&$2345$& $0.00192$&$14153$\\
Recursive 1& ${\bf 0.00318}$& $ 192$& $ {\bf 0.00215}$&$1209$&$ {\bf 0.00180}$&${\bf 7121}$\\
Recursive 2& $0.00334$& ${\bf 184}$& $0.00224$& ${\bf 1199}$&$ 0.00194$&$7144$\\
\hline
$\rho=80\%$&MWISE &CPU& MWISE &CPU& MWISE &CPU\\
\hline
Nonrecursive&$0.00334$&$366$& $0.00244$&$2543$& $0.00194$&$14532$\\
Recursive 1& ${\bf 0.00324}$&$194$& $ {\bf 0.00234}$&$1321$&$ {\bf 0.00184}$&$7342$\\
Recursive 2& $0.00338$&${\bf 188}$& $0.00246$&${\bf 1306}$&$ 0.00194$&${\bf 7224}$\\
\hline
& \multicolumn{6}{c}{$Y\sim \mathcal{W}eibull\left(2,1\right)$}\\
\hline
$\rho=30\%$&MWISE &CPU& MWISE &CPU& MWISE &CPU\\
\hline
Nonrecursive&$0.02342$&$345$& $0.02082$&$2442$& $0.01742$&$14121$\\
Recursive 1& ${\bf 0.02265}$&$173$& $ {\bf 0.001944}$&$1235$&$ {\bf 0.01638}$&$7214$\\
Recursive 2& $0.02348$&${\bf 163}$& $0.02105$&${\bf 1226}$&$ 0.01744$&${\bf 7116}$\\
\hline
$\rho=50\%$&MWISE &CPU& MWISE &CPU& MWISE &CPU\\
\hline
Nonrecursive&$0.02354$&$356$& $0.02142$&$2448$& $0.01754$&$14098$\\
Recursive 1& ${\bf 0.02289}$&${\bf 174}$& $ {\bf 0.02008}$&$1212$&$ {\bf 0.01644}$&${\bf 7032}$\\
Recursive 2& $ 0.02358$&$ 186$& $0.02144$&${\bf 1176}$&$ 0.01760$&$7034$\\
\hline
$\rho=80\%$&MWISE &CPU& MWISE &CPU& MWISE &CPU\\
\hline
Nonrecursive&$0.02372$&$352$& $0.02152$&$2356$& $0.001758$ &$14024$\\
Recursive 1& ${\bf 0.02314}$&${\bf 176}$& $ {\bf 0.02089}$&$1209$&$ {\bf 0.01654}$&$7014$\\
Recursive 2& $0.02369$&$ 184$& $0.02154$&${\bf 1165}$&$ 0.01760$&${\bf 7006}$\\
\hline
\end{tabular}
\end{center}
\caption{Quantitative comparison between the nonrecursive Horvitz-Thompson-type KDE~(\ref{eq:Nad}) and two recursive estimators; recursive $1$ correspond to the proposed estimator~(\ref{eq:fn}) with the choice of $\left(\gamma_n\right)=\left(n^{-1}\right)$, resursive $2$ correspond to the estimator~(\ref{eq:fn}) with the choice of $\left(\gamma_n\right)=\left(\left[4/5\right]n^{-1}\right)$.  We consider the proportion of missing at random equal to $70\%$ and the auxiliary variable $\rho$, equal to $30\%$ (in the first block), to $50\%$ (in the second block) and equal to $70\%$ (in the third block) of each considered distribution, we consider three sample sizes $n=100$, $n=200$ and $n=500$, the number of simulations is $500$, and we compute the Mean Weighted Integrated Squared Error  ($MWISE$) and the \texttt{CPU} time in seconds.}\label{Tab:4}
\end{table}

\begin{center}
\begin{figure}[!h]
    \includegraphics[width=0.65\textwidth,angle=270,clip=true,trim=40 0 0 0]{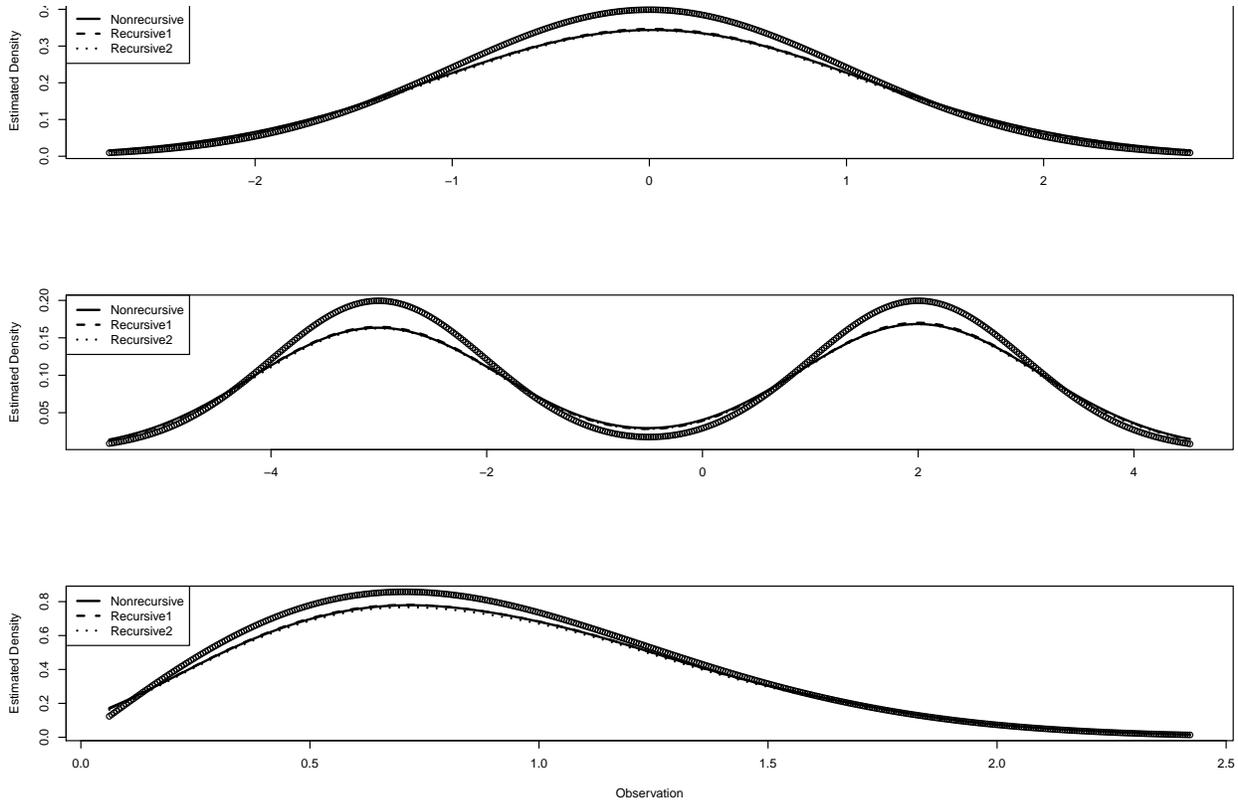}
\caption{Quantitative comparison between the nonrecursive Horvitz-Thompson-type KDE~(\ref{eq:Nad}) and two recursive estimators; recursive $1$ correspond to the proposed estimator~(\ref{eq:fn}) with the choice of $\left(\gamma_n\right)=\left(n^{-1}\right)$, resursive $2$ correspond to the estimator~(\ref{eq:fn}) with the choice of $\left(\gamma_n\right)=\left(\left[4/5\right]n^{-1}\right)$ for (a) normal, (b) normal mixture and (c) weibull distributions with $70\%$ missing data and correlation between the response and auxiliary variable of $\rho=0.8$. Based on $500$ simulated datasets of size $200$.}
\label{Fig:4}
\end{figure}
\end{center}

From Figures~\ref{Fig:1}~\ref{Fig:2},~\ref{Fig:3},~\ref{Fig:4}, Tables~\ref{Tab:1},~\ref{Tab:2},~\ref{Tab:3} and~\ref{Tab:4}, we conclude that
\begin{itemize}
\item The proposed recursive kernel density estimator~(\ref{eq:fn}), with the stepsize $\left(\gamma_n\right)=\left(n^{-1}\right)$ is closer to the true
density function as compared with the nonrecursive Horvitz-Thompson-type KDE~(\ref{eq:Nad}).
\item In all the considered densities, the proposed recursive KDE with the stepsize $\left(\gamma_n\right)=\left(n^{-1}\right)$ outperformed the nonrecursive Horvitz-Thompson-type KDE~(\ref{eq:Nad}) in terms of estimation error and computational costs.
 
\item The estimators get closer to the true density function as sample size increase and the proportion of the missing at random decrease.
\item The \texttt{CPU} time using the proposed recursive estimators are approximately two times faster than the \texttt{CPU} time using the nonrecursive Horvitz-Thompson-type KDE.
\end{itemize}

\subsubsection{Local density estimation}
In order to investigate the comparison between the proposed estimators in the case of local density estimation, we followed the paper of~\citet{Dut14} and we consider two sample sizes: $n=100$ and $1000$, the number of simulations is $N=200$. Moreover, we consider five densities functions $f$ : 1- the normal mixture distribution : $X\sim 1/2\mathcal{N}\left(-1,1/\sqrt{2}\right)+1/2\mathcal{N}\left(1,1/\sqrt{2}\right)$, 2- the normal mixture distribution : $X\sim 1/2\mathcal{N}\left(-1,1\right)+1/2\mathcal{N}\left(1,1\right)$, 3- the standard normal : $X\sim \mathcal{N}\left(0,1\right)$, 4- the exponetial distribution with rate $1$: $X\sim \mathcal{E}xp\left(1\right)$, 5- the Cauchy distribution with location $0$ and scale $1$. We report in Table~\ref{Tab:0}, in the first two columns the bandwidth and the Monte Carlo estimates of the ratio of the square root of the MSE of the estimator proposed by~\citet{Cha10} (see Table 1 of~\citet{Dut14}), then, we give the bandwidth and the Monte Carlo estimates of the ratio of the square root of the MSE of the nonrecursive Horvitz-Thompson-type KDE using the proposed plug-in bandwidth selection~(\ref{hoptim:rose:plug:in:loc}) and in the last two columns we give the Monte Carlo estimates of the square ratio of the root of the MSE of the proposed recursive estimator using the proposed plug-in bandwidth selection~(\ref{hoptim:gamma:loc:plug:in}). 
\begin{table}
\begin{center}
\begin{tabular}{llcccccc}
\hline
Density&$\left(x_0,n\right)$ &$h_{CLP}$ &CLP &$h_{NRPI}$& NRPI & $h_{RPI}$ &RPI\\
\hline
$\frac{1}{2}\mathcal{N}\left(-1,\frac{1}{\sqrt{2}}\right)+\frac{1}{2}\mathcal{N}\left(1,\frac{1}{\sqrt{2}}\right)$&$\left(0,100\right)$& $0.192$&$0.227$&$0.606$& $0.191$& $0.716$& $0.189$\\
&$\left(0,1000\right)$& $0.151$& $0.233$&$0.274$&$0.086$&$0.579$&$0.158$\\
$\frac{1}{2}\mathcal{N}\left(-1,1\right)+\frac{1}{2}\mathcal{N}\left(1,1\right)$&$\left(0,100\right)$& $0.099$&$0.118$&$0.806$&$0.083$&$0.865$&$0.085$\\
&$\left(0,1000\right)$& $0.039$ &$0.048$&$0.474$&$0.037$&$0.751$&$0.046$\\
$\mathcal{N}\left(0,1\right)$&$\left(-1.282,100\right)$& $0.149$& $0.175$&$0.558$&$0.161$&$0.757$&$0.135$\\
&$\left(-1.282,1000\right)$&$0.070$& $0.131$&$0.321$&$0.071$&$0.445$&$0.067$\\
&$\left(0,100\right)$& $0.131$&$0.107$&$0.408$ &$0.170$ &$0.646$ &$0.133$\\
&$\left(0,1000\right)$&$0.056$&$0.045$ &$0.334$ &$0.069$ &$0.560$ &$0.074$\\
&$\left(1.282,100\right)$& $0.168$&$0.178$&$0.566$ &$0.148$ &$0.701$ &$0.134$\\
&$\left(1.282,1000\right)$& $0.070$&$0.129$&$0.353$&$0.069$&$0.523$&$0.069$\\
$\mathcal{E}xp\left(1\right)$&$\left(0.1054,100\right)$& $0.435$&$0.180$&$0.177$&$0.338$&$0.234$&$0.403$\\
&$\left(0.1054,1000\right)$& $0.260$&$0.088$&$0.083$&$0.130$&$0.196$&$0.357$\\
&$\left(0.693,100\right)$& $0.145$&$0.135$&$0.387$&$0.071$&$0.527$&$0.118$\\
&$\left(0.693,1000\right)$& $0.070$&$0.071$&$0.301$&$0.038$&$0.411$&$0.029$\\
&$\left(2.303,100\right)$& $0.213$&$0.232$&$0.406$&$0.267$&$0.300$&$0.296$\\
&$\left(2.303,1000\right)$& $0.077$&$0.101$&$0.277$&$0.101$&$0.309$&$0.100$\\
$Cauchy\left(0,1\right)$&$\left(-3.078,100\right)$& $0.390$& $1.573$&$0.350$&$0.533$&$0.471$&$0.450$\\
&$\left(-3.078,1000\right)$&$0.147$& $1.514$&$0.205$& $0.203$& $0.358$& $0.150$\\
&$\left(0,100\right)$& $0.218$&$0.135$&$0.332$& $0.152$& $0.447$& $0.163$\\
&$\left(0,1000\right)$& $0.076$&$0.046$&$0.190$& $0.070$& $0.305$& $0.088$\\
&$\left(3.078,100\right)$&$0.435$& $1.662$&$0.336$&$0.456$&$0.398$&$0.415$\\
&$\left(3.078,1000\right)$& $0.151$&$1.513$&$0.213$&$0.211$&$0.471$&$0.149$\\
\hline
\end{tabular}
\end{center}
\caption{Monte Carlo estimates of $\frac{\sqrt{MSE\left(f_n\left(x_0\right)\right)}}{f\left(x_0\right)}$ and the bandwidth obtained respectively by CLP method, the nonrecursive Horvitz-Thompson-type KDE~(\ref{eq:Nad}) using the proposed plug-in method and the proposed recursive estimator using the proposed plug-in method.}\label{Tab:0}
\end{table}
In Table~\ref{Tab:0}, we consider the estimation of (a) the $1/2\mathcal{N}\left(-1,1/\sqrt{2}\right)+1/2\mathcal{N}\left(1,1/\sqrt{2}\right)$ density at $x_0=0$, (b) the $1/2\mathcal{N}\left(-1,1\right)+1/2\mathcal{N}\left(1,1\right)$ density at $x_0=0$, and the estimation of (c) $\mathcal{N}\left(0,1\right)$, (d) $\mathcal{E}xp\left(1\right)$, and (e) Cauchy$\left(0,1\right)$ densities at $x_0$ equal to the $10$th, $50$th, and $90$th percentiles.\\
From Table~\ref{Tab:0}, we have the following observations.
\begin{description}
\item[(i)] The proposed recursive estimator and the Horvitz-Thompson-type KDE, preformed better than the CLP method in many situations: using the mixed normal densities ((a) and (b)), using the standard normal density (c) and standard Cauchy (e) in the tail regions ($10$th and $90$th percentiles), using the exponential density (d) with $x_0$ equal to the median. 
\item[(ii)] The proposed recursive estimator and the Horvitz-Thompson-type KDE, seems to be consistent in all the cases (its MISE seems to decrease as $n$ increased), however the CLP estimator does not seem to be consistent (see for example the mixed normal density (a)).
\item[(ii)] No one between the proposed recursive estimator and the Horvitz-Thompson-type KDE can be claimed to be the best in all the considered cases. 
\end{description}

\subsection{Examples}\label{subsection:real}
In this section we report our analysis of the two datasets mentioned in the Introduction, the first one concerne the simulated data of the Aquitaine cohort of $\texttt{HIV}$-1 infected patients and the second one concerne the coriell cell lines using the chromosome number 11.
\paragraph{Simulated data : Aquitaine cohort of $\texttt{HIV}$-1 infected patients}

The aim of this part was to compare the performance of the proposed recursive KDE to the nonrecursive Horvitz-Thompson-type KDE using the Aquitaine cohort of $\texttt{HIV}$-1 infected patients, see the paper~\citet{Thi00} and receiving highly active anti-retroviral therapy.

We generated $200$ samples of $100$ subjects. The number of repeated measures for each subject was randomly distributed between $2$ and $7$ (mean $4$) and the times of measurements were uniformly distributed between $0$ and $6$.

For comparisons, we simulated data according to the following linear mixed effects model:  
\begin{eqnarray*}
Y_{i,k}=\texttt{A1}+\texttt{A2}*X_{i,k}+b_{i}+\varepsilon_{i,k}
\end{eqnarray*}
with the assumption that $b_{i}\sim\mathcal{N}\left(0,\sigma_1^2\right)$ and $\varepsilon_{i,k}\sim\mathcal{N}\left(0,\sigma^2\right)$. We assumed that random coefficient $\left(b_{i}, \varepsilon_{i,k}\right)$ were independent of each other. The values of parameters were chosen to be similar to those obtained from the Aquitaine cohort of $\texttt{HIV}$-1 infected patients.  The parameters were $\texttt{A1}=4$, $\texttt{A2}=-0.5$, $\sigma_1^2=0.25$ and $\sigma^2=1$.

Then, we fixed the proportion of missing at random to be equal to $30\%$. 
    
\begin{center}
\begin{figure}[!h]
    \includegraphics[width=0.65\textwidth,angle=270,clip=true,trim=40 0 0 0]{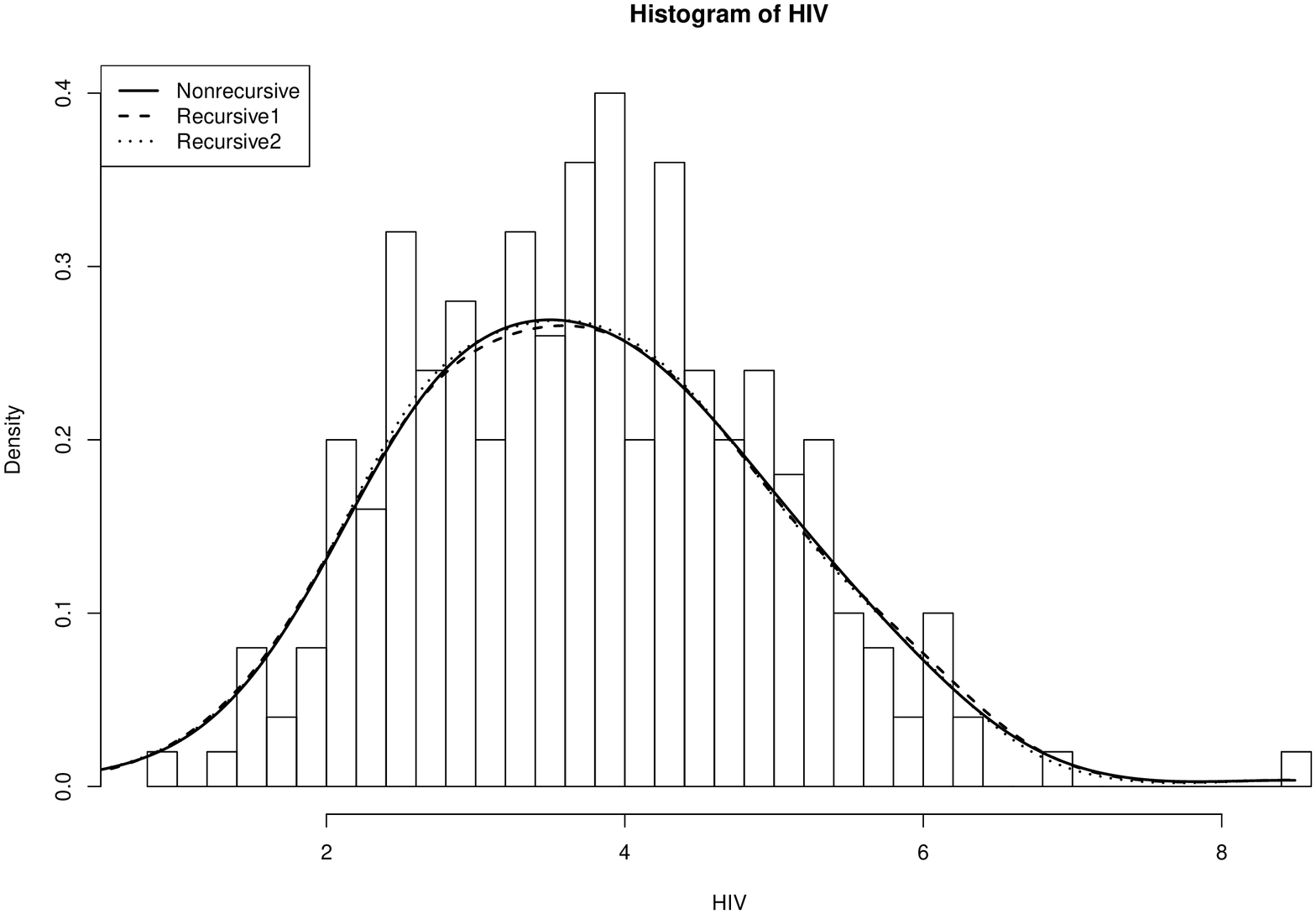}
\caption{Quantitative comparison between the nonrecursive Horvitz-Thompson-type KDE~(\ref{eq:Nad}) and two recursive estimators; recursive $1$ correspond to the proposed estimator~(\ref{eq:fn}) with the choice of $\left(\gamma_n\right)=\left(n^{-1}\right)$, resursive $2$ correspond to the estimator~(\ref{eq:fn}) with the choice of $\left(\gamma_n\right)=\left(\left[4/5\right]n^{-1}\right)$. Here we consider the fibroblast cell strains using the chromosome number 11 with the proportion of missing at random equal to $30\%$. Running time was roughly $483$s using the Horvitz-Thompson-type KDE, $248$s using recursive $1$ and $253$s using recursive $2$}
\label{Fig:5}
\end{figure}
\end{center}
In Figure~\ref{Fig:5}, we qualitatively compared the proposed recursive KDE using the two discussed stepsizes to the nonrecursive Horvitz-Thompson-type KDE. We observed that even when $30\%$ of the original measurements are missing, the proposed recursive estimators remain very accurate thus demontrating the effectiveness of our approach.

Moreover, our proposed recursive KDE was over two times faster than the nonrecursive Horvitz-Thompson-type KDE algorithm~(\ref{eq:Nad}).\\
Our second application concerned real dataset, for which we compared the proposed recursive KDE approach to the nonrecursive Horvitz-Thompson-type KDE algorithm using a special data using essentially in the context of change point problem.\\ 
\paragraph{Real dataset: Coriell cell lines using the chromosome number 11}

The data correspond to two array Comparative Genomic Hybridization (CGH) studies of Coriell cell lines, which appears in the \texttt{DNAcopy} package and \texttt{bcp} package (see the paper~\citet{Erd07}), and for more details about the data we advise the reader to see the paper of~\citet{Sni01}.

The data correspond to two array CGH studies of fibroblast cell strains. In particular in~\citet{Ven07}, they chose the studies {\bf GM05296} and {\bf GM13330}. After selecting only the mapped data from chromosomes $1$-$22$ and $X$, there are $2271$ data points. Here we perform an analysis on the {\bf GM13330} array CGH study described above. We simulated three different proportions of missing at random.

\begin{figure}
\begin{center}
\includegraphics[width=0.65\textwidth,angle=270,clip=true,trim=40 0 0 0]{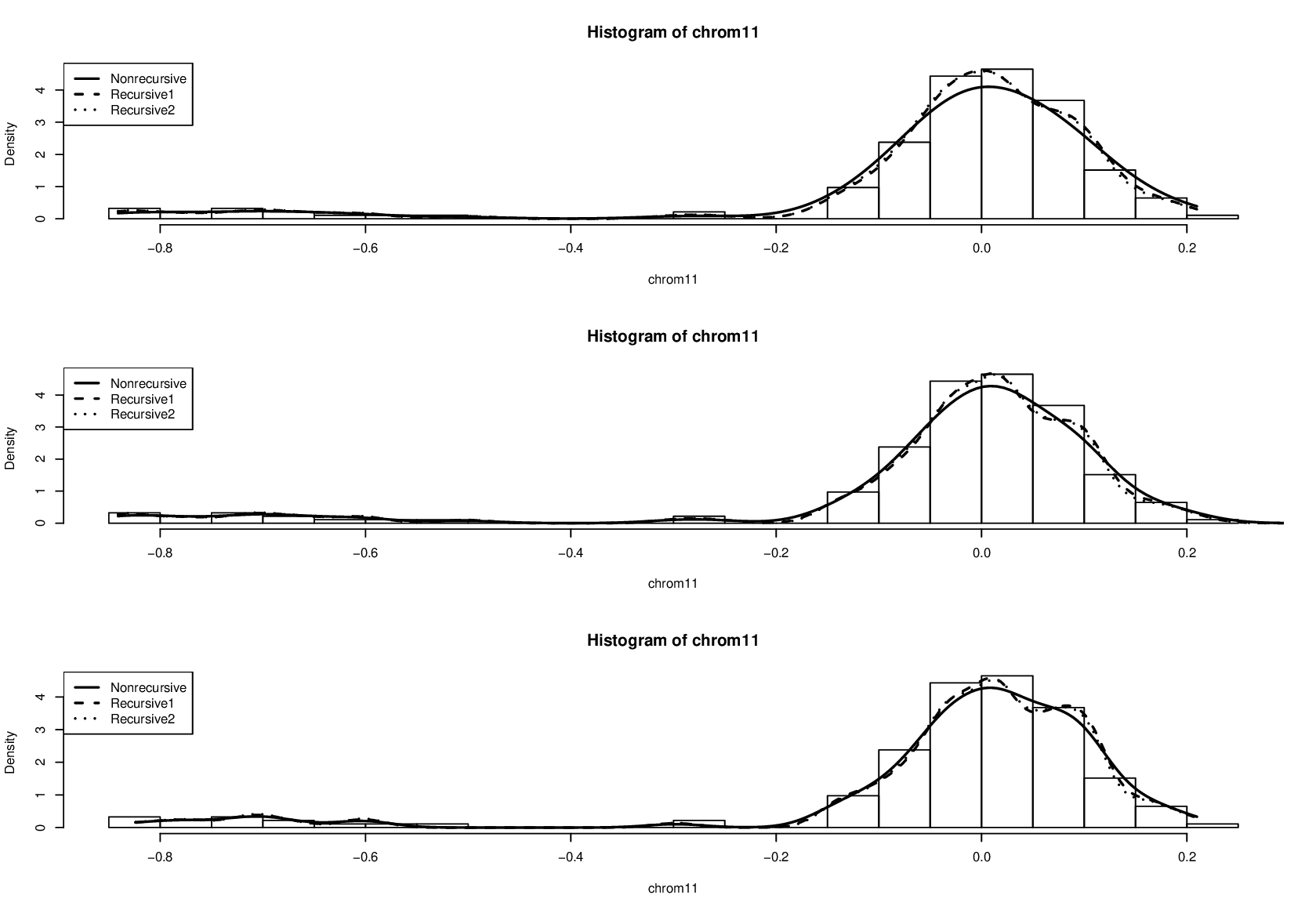}
\end{center}
\caption{Quantitative comparison between the nonrecursive Horvitz-Thompson-type KDE~(\ref{eq:Nad}) and two recursive estimators; recursive $1$ correspond to the proposed estimator~(\ref{eq:fn}) with the choice of $\left(\gamma_n\right)=\left(n^{-1}\right)$, resursive $2$ correspond to the estimator~(\ref{eq:fn}) with the choice of $\left(\gamma_n\right)=\left(\left[4/5\right]n^{-1}\right)$. Here we consider the fibroblast cell strains using the chromosome number 11 with the proportion of missing at random equal respectively to $0\%$ (in the first block), equal to $15\%$ (in the second block), and equal to $30\%$ (in the last block). Running time was roughly $211$s using the Horvitz-Thompson-type KDE, $107$s using recursive $1$ and $113$s using recursive $2$, when the proportion of missing at random equal to $30\%$}
\label{Fig:6}
\end{figure}
It is interesting to note from Figure~\ref{Fig:6} that the three estimators are quite similar, and closely follow the shape of the histogram. Moreover, we can claim that the proposed recursive estimators with the choice of the stepsize $\left(\gamma_n\right)=\left(n^{-1}\right)$ is closer to the shape of the histogram of the true data as compared with the nonrecursive Horvitz-Thompson-type KDE~(\ref{eq:Nad}).

\section{Conclusion}\label{section:conclusion}

This paper propose an automatic bandwidth selection of the recursive KDE under missing data~(\ref{eq:fn}). The proposed estimators are asymptotically follows normal distribution. The proposed estimators are compared to the nonrecursive Horvitz-Thompson-type KDE algorithm under missing data~(\ref{eq:Nad}). We showed that, using a recursive version of the Nadaraya-Watson estimator to estimate the propensity scores and using a specific data-driven bandwidth and some particularly stepsizes, the proposed recursive estimators can give better results compared to the nonrecursive Horvitz-Thompson-type estimator in terms of estimation error in the case of global estimation and quite similar results in the case of local estimation. Noting that the two estimators ((\ref{eq:fn}) and (\ref{eq:Nad})) can deal effectively with both missing and complete data.
The simulation study confirms the nice features of our proposed recursive estimators and satisfactory improvement in the \texttt{CPU} time in comparison to the nonrecursive Horvitz-Thompson-type KDE estimator.

In conclusion, the proposed recursive estimators allowed us to obtain quite better results compared to the nonrecursive Horvitz-Thompson-type KDE under missing data in terms of estimation error and much better in terms of computational costs. 
\appendix
\section{Technical proofs for the recursive KDE} \label{section:proofs}
Throughout this section we use the following notation:
\begin{eqnarray}
Q_n&=&\prod_{j=1}^{n}\left(1-\gamma_j\right)\nonumber,\\
Z_n\left(x\right)&=&h_n^{-1}\delta_n\pi_n^{-1}K\left(\frac{x-X_n}{h_n}\right).\label{eq:Zn}
\end{eqnarray}
Let us first state the following technical lemma. 

\begin{lemma}\label{L:1} 
Let $\left(v_n\right)\in \mathcal{GS}\left(v^*\right)$, 
$\left(\eta_n\right)\in \mathcal{GS}\left(-\eta\right)$, and $m>0$ 
such that $m-v^*\xi>0$ where $\xi$ is defined in~\eqref{eq:xi}. We have 
\begin{eqnarray*}
\lim_{n \to +\infty}v_nQ_n^{m}\sum_{k=1}^nQ_k^{-m}\gamma_kv_k^{-1}
=\left(m-v^*\xi\right)^{-1}. 
\end{eqnarray*}
Moreover, for all positive sequence $\left(\alpha_n\right)$ such that $\lim_{n \to +\infty}\alpha_n=0$, and all $C \in \mathbb{R}$,
\begin{eqnarray*}
\lim_{n \to +\infty}v_nQ_n^{m}\left[\sum_{k=1}^n Q_k^{-m} \eta_kv_k^{-1}\alpha_k+C\right]=0.
\end{eqnarray*}
\end{lemma}

Lemma \ref{L:1} is widely applied throughout the proofs. Let us underline that it is its application, which requires Assumption $(A2)(iii)$ on the limit of $(n\gamma_n)$ as $n$ goes to infinity.\\

Our proofs are organized as follows. Propositions \ref{prop:bias:var:lambn} and \ref{prop:MISE:lamb} in Sections \ref{Section:prop:bias:var:lambn} 
and \ref{Section:prop:MISE:lambn} respectively, Theorem \ref{theo:TLC1} in Section \ref{Section:theo:TLC}.

\subsection{Proof of Proposition~\ref{prop:bias:var:lambn}} \label{Section:prop:bias:var:lambn}

In view of~(\ref{eq:fn}) and~(\ref{eq:Zn}), we have
\begin{eqnarray*}
f_n\left(x\right)-f\left(x\right)
&=&\left(1-\gamma_n\right)\left(f_{n-1}\left(x\right)-f\left(x\right)\right)+\gamma_n\left(Z_n\left(x\right)-f\left(x\right)\right)\nonumber\\
&=&\sum_{k=1}^{n-1}\left[\prod_{j=k+1}^n\left(1-\gamma_j\right)\right]\gamma_k\left(Z_k\left(x\right)-f\left(x\right)\right)+\gamma_n\left(Z_n\left(x\right)-f\left(x\right)\right)
\nonumber\\
&&+\left[\prod_{j=1}^n\left(1-\gamma_j\right)\right]\left(f_0\left(x\right)-f\left(x\right)\right)\nonumber\\
&=&Q_n\sum_{k=1}^nQ_k^{-1}\gamma_k\left(Z_k\left(x\right)-f\left(x\right)\right)+Q_n\left(f_0\left(x\right)-f\left(x\right)\right).
\end{eqnarray*} 
It follows that
\begin{eqnarray*}
\mathbb{E}\left(f_{n}\left(x\right)\right)-f\left(x\right)
&=&Q_n\sum_{k=1}^nQ_k^{-1}\gamma_k\left(\mathbb{E}\left(Z_k\left(x\right)\right)-f\left(x\right)\right)+Q_n\left(f_0\left(x\right)-f\left(x\right)\right).
\end{eqnarray*} 
Moreover, we have
\begin{eqnarray}\label{eq:Znp}
\mathbb{E}\left[Z_k^p\left(x\right)\right]&=&h_k^{-p}\pi_k^{-p}\mathbb{E}\left[\delta_kK^p\left(\frac{x-X_k}{h_k}\right)\right]\nonumber\\
&=&h_k^{-p}\pi_k^{-p}\mathbb{E}\left[\mathds{1}_{\left\{T_k=X_k\right\}}K^p\left(\frac{x-X_k}{h_k}\right)\right]\nonumber\\
&=&h_k^{-p}\pi_k^{-p}\mathbb{P}\left[T_k= X_k\right]\mathbb{E}\left[K^p\left(\frac{x-X_k}{h_k}\right)\vert \left\{T_k=X_k\right\}\right]\nonumber\\
&=&h_k^{-p}\pi_k^{-p}\mathbb{P}\left[\delta_k=1\vert T_i\right] \mathbb{E}\left[K^p\left(\frac{x-T_k}{h_k}\right)\right]\nonumber\\
&=&h_k^{-p}\pi_k^{-p+1}\mathbb{E}\left[K^p\left(\frac{x-T_k}{h_k}\right)\right]\nonumber\\
&=&h_k^{-p}\pi_k^{-p+1}\int_{\mathbb{R}}K^p\left(\frac{x-t}{h_k}\right)f\left(t\right)dt\nonumber\\
&=&h_k^{-p+1}\pi_k^{-p+1}\int_{\mathbb{R}}K^p\left(z\right)f\left(x-zh_k\right)dz
\end{eqnarray}
Then, it follows from~(\ref{eq:Znp}), for $p=1$, that
\begin{eqnarray}
\mathbb{E}\left[Z_k\left(x\right)\right]-f\left(x\right)&=&\int_{\mathbb{R}}K\left(z\right)\left[f\left(x-zh_k\right)-f\left(x\right)\right]dz\nonumber\\
&=&\frac{h_k^2}{2}f^{\left(2\right)}\left(x\right)\mu_2\left(K\right)+\eta_k\left(x\right),\label{eq:EZn}
\end{eqnarray}
with
\begin{eqnarray*}
\eta_k\left(x\right)=\int_{\mathbb{R}}K\left(z\right)\left[f\left(x-zh_k\right)-f\left(x\right)-\frac{1}{2}z^2h_k^2f^{\left(2\right)}\left(x\right)\right]dz,
\end{eqnarray*}
and, since $f$ is bounded and continuous at $x$, we have $\lim_{k\to \infty}\eta_k\left(x\right)=0$. In the case $a\leq \alpha/5$, we have $\lim_{n\to \infty}\left(n\gamma_n\right)>2a$; the application of Lemma~\ref{L:1} then gives 
\begin{eqnarray*}
\mathbb{E}\left[f_n\left(x\right)\right]-f\left(x\right)&=&\frac{1}{2}f^{\left(2\right)}\left(x\right)\int_{\mathbb{R}}z^2K\left(z\right)dzQ_n\sum_{k=1}^nQ_k^{-1}\gamma_kh_k^2\left[1+o\left(1\right)\right]\\
&&+Q_n\left(f_0\left(x\right)-f\left(x\right)\right)\nonumber\\
&=&\frac{1}{2\left(1-2a\xi\right)}f^{\left(2\right)}\left(x\right)\mu_2\left(K\right)\left[h_n^2+o\left(1\right)\right],
\end{eqnarray*}
and~(\ref{bias:lambn}) follows. In the case $a>\alpha/5$, we have $h_n^2=o\left(\sqrt{\gamma_nh_n^{-1}}\right)$, and $\lim_{n\to \infty}\left(n\gamma_n\right)>\left(\alpha-a\right)/2$, then Lemma~\ref{L:1} ensures that
\begin{eqnarray*}
\mathbb{E}\left[f_n\left(x\right)\right]-f\left(x\right)&=&Q_n\sum_{k=1}^nQ_k^{-1}\gamma_ko\left(\sqrt{\gamma_kh_k^{-1}}\right)+O\left(Q_n\right)\\
&=&o\left(\sqrt{\gamma_nh_n^{-1}}\right).
\end{eqnarray*}
which gives~(\ref{bias:lambn:bis}). 
Further, we have 
\begin{eqnarray}\label{eq:varzk}
Var\left[f_n\left(x\right)\right]&=&Q_n^2\sum_{k=1}^nQ_k^{-2}\gamma_k^2Var\left[Z_k\left(x\right)\right]\nonumber\\
&=&Q_n^2\sum_{k=1}^nQ_k^{-2}\gamma_k^2\left(\mathbb{E}\left(Z_k^2\left(x\right)\right)-\left(\mathbb{E}\left(Z_k\left(x\right)\right)\right)^2\right).
\end{eqnarray}
Moreover, in view of~(\ref{eq:Znp}), for $p=2$, that 
\begin{eqnarray}
\mathbb{E}\left(Z_k^2\left(x\right)\right)&=&h_k^{-1}\pi_k^{-1}\int_{\mathbb{R}}f\left(x-zh_k\right)K^2\left(z\right)dz\nonumber\\
&=&h_k^{-1}\pi_k^{-1}f\left(x\right)\int_{\mathbb{R}}K^2\left(z\right)dz+\nu_k\left(x\right),\label{eq:EZk2}
\end{eqnarray}
with 
\begin{eqnarray*}
\nu_k\left(x\right)=h_k^{-1}\pi_k^{-1}\int_{\mathbb{R}}K^2\left(z\right)\left[f\left(x-zh_k\right)-f\left(x\right)\right]dz.
\end{eqnarray*}
Moreover, it follows from~(\ref{eq:EZn}), that
\begin{eqnarray}
\mathbb{E}\left[Z_k\left(x\right)\right]&=&f\left(x\right)+\widetilde{\nu}_k\left(x\right),\label{eq:EZk1}
\end{eqnarray}
with
\begin{eqnarray*}
\widetilde{\nu}_k\left(x\right)=\int_{\mathbb{R}}K\left(z\right)\left[f\left(x-zh_k\right)-f\left(x\right)\right]dz.
\end{eqnarray*}
Then, it follows from~(\ref{eq:varzk}), (\ref{eq:EZk2}) and (\ref{eq:EZk1}), that
\begin{eqnarray*}
Var\left[f_n\left(x\right)\right]&=&f\left(x\right)R\left(K\right)Q_n^2\sum_{k=1}^nQ_k^{-2}\gamma_k^2h_k^{-1}\pi_k^{-1}+Q_n^2\sum_{k=1}^nQ_k^{-2}\gamma_k^2\nu_k\left(x\right)\nonumber\\
&&-f^2\left(x\right)Q_n^2\sum_{k=1}^nQ_k^{-2}\gamma_k^2-2f\left(x\right)Q_n^2\sum_{k=1}^nQ_k^{-2}\gamma_k^2\widetilde{\nu}_k\left(x\right)\nonumber\\
&&-Q_n^2\sum_{k=1}^nQ_k^{-2}\gamma_k^2\widetilde{\nu}_k^2\left(x\right).
\end{eqnarray*}
Since $f$ is bounded continuous, we have $\lim_{k\to \infty}\nu_k\left(x\right)=0$ and $\lim_{k\to \infty}\widetilde{\nu_k}\left(x\right)=0$. In the case $a\geq \alpha/5$, we have $\lim_{n\to \infty}\left(n\gamma_n\right)>\left(\alpha-a\right)/2$, and the application of Lemma~\ref{L:1} gives
\begin{eqnarray*}
Var\left[f_n\left(x\right)\right]&=&\gamma_nh_n^{-1}\pi_n^{-1}\left(2-\left(\alpha-a\right)\xi\right)^{-1}f\left(x\right)R\left(K\right)
+o\left(\gamma_nh_n^{-1}\right),
\end{eqnarray*}
which proves~(\ref{var:lambn}). Now, in the case $a<\alpha/5$, we have $\gamma_nh_n^{-1}=o\left(h_n^4\right)$, and $\lim_{n\to \infty}\left(n\gamma_n\right)>2a$, then the application of Lemma~\ref{L:1} gives 
\begin{eqnarray*}
Var\left[f_n\left(x\right)\right]&=&Q_n^2\sum_{k=1}^nQ_k^{-2}\gamma_ko\left(h_k^4\right)\\
&=&o\left(h_n^4\right),
\end{eqnarray*}
which proves~(\ref{var:lambn:bis}).
\subsection{Proof of Proposition~\ref{prop:MISE:lamb}} \label{Section:prop:MISE:lambn}
Following similar steps as the proof of the Proposition 2 of~\citet{Mok09}, we proof the Proposition~\ref{prop:MISE:lamb}.

\subsection{Proof of Theorem~\ref{theo:TLC1}} \label{Section:theo:TLC}
Let us at first assume that, if $a\geq\alpha/5$, then 
\begin{eqnarray}\label{eq:22}
\sqrt{\gamma_n^{-1} h_n\pi_n}\left(f_{n}\left( x\right)-\mathbb{E}\left[f_n\left(x\right)\right]\right)\stackrel{\mathcal{D}}{\rightarrow}\mathcal{N}\left( 0,
\left(2-\left(\alpha-a\right)\xi\right)^{-1}f\left(x\right)R\left(K\right)\right).
\end{eqnarray}
In the case when $a>\alpha/5$, Part 1 of Theorem~\ref{theo:TLC1} follows from the combination of~\eqref{bias:lambn:bis} and~\eqref{eq:22}. In the case when $a=\alpha/5$, Parts 1 and 2 of Theorem~\ref{theo:TLC1} follow from the combination of~\eqref{bias:lambn} and~\eqref{eq:22}. In the case $a<\alpha/5$,~\eqref{var:lambn:bis} implies that 
\begin{eqnarray*}
h_n^{-2}\left(f_n\left(x\right)-\mathbb{E}\left(f_n\left(x\right)\right)\right)\stackrel{\mathbb{P}}{\rightarrow}0,
\end{eqnarray*}
and the application of~\eqref{bias:lambn} gives Part 2 of Theorem~\ref{theo:TLC1}.\\

We now prove~\eqref{eq:22}. 
In view of~\eqref{eq:fn}, we have
\begin{eqnarray*}
f_n\left(x\right)-\mathbb{E}\left[f_n\left(x\right)\right]=Q_n\sum_{k=1}^nQ_k^{-1}\gamma_k\left(Z_k\left(x\right)-\mathbb{E}\left[Z_k\left(x\right)\right]\right).
\end{eqnarray*}
Set
\begin{eqnarray*}
Y_k\left(x\right)&=&\Pi_k^{-1}\gamma_k\left(Z_k\left(x\right)-\mathbb{E}\left[Z_k\left(x\right)\right]\right).
\end{eqnarray*}
The application of Lemma~\ref{L:1} ensures that
\begin{eqnarray*} 
v_n^2&=&\sum_{k=1}^nVar\left(Y_{k}\left(x\right)\right)\nonumber\\
&=&\sum_{k=1}^n\Pi_k^{-2}\gamma_k^2Var\left(Z_k\left(x\right)\right)\nonumber\\
&=&\sum_{k=1}^n\Pi_k^{-2}\gamma_k^2h_k^{-1}\pi_k^{-1}\left[f\left(x\right)R\left(K\right)+o\left(1\right)\right]\nonumber\\
&=&Q_n^{-2}\gamma_nh_n^{-1}\pi_n^{-1}\left[\left(2-\left(\alpha-a\right)\xi\right)^{-1}f\left(x\right)R\left(K\right)+o\left(1\right)\right].
\end{eqnarray*}
On the other hand, we have, for all $p>0$, 
\begin{eqnarray*}
\mathbb{E}\left[\left|Z_k\left(x\right)\right|^{2+p}\right] &=&
O\left(h_k^{-1-p}\right),
\end{eqnarray*}
and, since $\lim_{n\to\infty}\left(n\gamma_n\right)>\left(\alpha-a\right)/2$, there exists $p>0$ such that $\lim_{n\to \infty}\left(n\gamma_n\right)>\frac{1+p}{2+p}\left(\alpha-a\right)$. Applying Lemma $\ref{L:1}$, we get 
\begin{eqnarray*}
\sum_{k=1}^n\mathbb{E}\left[\left|Y_{k}\left(x\right)\right|^{2+p}\right]&=&O\left(\sum_{k=1}^n Q_k^{-2-p}\gamma_k^{2+p}\mathbb{E}\left[\left|Y_k\left(x\right)\right|^{2+p}\right]\right)\nonumber\\
&=&O\left(\sum_{k=1}^n \frac{\Pi_k^{-2-p}\gamma_k^{2+p}}{h_k^{1+p}}\right)\\
&=&O\left(\frac{\gamma_n^{1+p}}{Q_n^{2+p}h_n^{1+p}}\right)\nonumber,
\end{eqnarray*}
and we thus obtain 
\begin{eqnarray*}
\frac{1}{v_n^{2+p}}\sum_{k=1}^n\mathbb{E}\left[\left|Y_{k}\left(x\right)\right|^{2+p}\right]& = & O\left({\left[\gamma_nh_n^{-1}\right]}^{p/2}\right)=o\left(1\right).
\end{eqnarray*}
The convergence in~\eqref{eq:22} then follows from the application of Lyapounov's Theorem.

\section{Results for Horvitz-Thompson-type KDE}\label{section:results:HTKDE}
\subsection{Global density estimation using HT KDE}\label{sub:sec:glo:HTKDE}
Now, let us recall that the bias and variance of the nonrecursive Horvitz-Thompson-type KDE defined by~(\ref{eq:Nad}) are given by
\begin{eqnarray*}
\mathbb{E}\left[\widetilde{f}_n\left(x\right)\right]-f\left(x\right)
=\frac{h_n^2}{2}f^{\left(2\right)}\left(x\right)\mu_2\left(K\right)+o\left(h_n^2\right),
\end{eqnarray*}
and
\begin{eqnarray*}
Var\left[\widetilde{f}_n\left(x\right)\right]&=&\frac{\pi_n^{-1}}{nh_n}f\left(x\right)R\left(K\right)+o\left(\frac{1}{nh_n}\right).
\end{eqnarray*}
It follows that,
\begin{eqnarray*} AMWISE\left[\widetilde{f}_n\right]&=&\frac{\pi_n^{-1}}{nh_n}I_1R\left(K\right)+\frac{h_n^4}{4}I_2h_n^4\mu_2^2\left(K\right).
\end{eqnarray*}
Then, to minimize the $AMWISE$ of $\widetilde{f}_n$, the bandwidth $\left(h_n\right)$ must equal to
\begin{eqnarray}\label{hoptim:rose}
\left(\left(\frac{I_1}{I_2}\right)^{1/5}\left\{\frac{R\left(K\right)}{\mu_2^2\left(K\right)}\right\}^{1/5}\pi_n^{-1/5}n^{-1/5}\right),
\end{eqnarray}
and we have
\begin{eqnarray*} AMWISE\left[\widetilde{f}_n\right]=\frac{5}{4}I_1^{4/5}I_2^{1/5}\Theta\left(K\right)\pi_n^{-4/5}n^{-4/5}.
\end{eqnarray*}
To estimate the optimal bandwidth~(\ref{hoptim:rose}), we must estimate $I_1$ and $I_2$. For this purpose, we use the following two kernel estimators :
\begin{eqnarray}
\widetilde{I}_1&=&\frac{1}{n\left(n-1\right)b_n}\sum_{\substack{i,j=1\\i\not=j}}^n\delta_j\widehat{\pi}_{NWj}^{-1}K_b\left(\frac{X_i-X_j}{b_n}\right),\label{I1:norec}\\
\widetilde{I}_2&=&\frac{1}{n^3b^{\prime 6}_n}\sum_{\substack{i,j,k=1\\j\not=k}}^n\delta_j\delta_k\widehat{\pi}_{NWj}^{-1}\widehat{\pi}_{NWk}^{-1}K_{b^{\prime}}^{\left(2\right)}\left(\frac{X_i-X_j}{b^{\prime}_n}\right)K_{b^{\prime}}^{\left(2\right)}\left(\frac{X_i-X_k}{b_n^{\prime}}\right)\label{I2:norec}.
\end{eqnarray}
where $K_b$ and $K_{b^{\prime}}$ are a kernels, $b_n$ and $b^{\prime}_n$ are respectively the associated bandwidth given in~(\ref{eq:h:initial}).\\ 
We showed that in order to minimize the $AMISE$ of $\widetilde{I}_1$ respectively of $\widetilde{I}_2$, the pilot bandwidth $\left(b_n\right)$ respectively $\left(b_n\right)$  must belong to $\mathcal{GS}\left(-2/5\right)$, respectively to $\mathcal{GS}\left(-3/14\right)$.\\
Then the plug-in estimator of the bandwidth $\left(h_n\right)$ using the nonrecursive estimator~(\ref{eq:Nad}), is given by
\begin{eqnarray}\label{hoptim:rose:plug:in}
\left(\left(\frac{\widetilde{I}_1}{\widetilde{I}_2}\right)^{1/5}\left\{\frac{R\left(K\right)}{\mu_2^2\left(K\right)}\right\}^{1/5}\widehat{\pi}_{NWn}^{-1/5}n^{-1/5}\right),
\end{eqnarray}
and the plug-in of the $AMWISE$ of the nonrecursive estimator~(\ref{eq:Nad}), is given by
\begin{eqnarray*}
\widetilde{AMWISE}\left[\widetilde{f}_n\right]=\frac{5}{4}\widetilde{I}_1^{4/5}\widetilde{I}_2^{1/5}\Theta\left(K\right)\widehat{\pi}_{NWn}^{-4/5}n^{-4/5}.
\end{eqnarray*}
\subsection{Local density estimation using the HT KDE}\label{sub:sec:loc:HTKDE}
Now, using the nonrecursive Horvitz-Thompson-type KDE and in order to minimize the $AMSE$ of $\widetilde{f}_n$, the bandwidth $\left(h_n\right)$ must equal to 
\begin{eqnarray}\label{hoptim:rose}
\left(\left(\frac{f\left(x\right)}{\left(f^{\left(2\right)}\left(x\right)\right)^2}\right)^{1/5}\left\{\frac{R\left(K\right)}{\mu_2^2\left(K\right)}\right\}^{1/5}\pi_n^{-1/5}n^{-1/5}\right).
\end{eqnarray} 
Moreover, since the $f\left(x\right)$ and $f^{\left(2\right)}\left(x\right)$ are not known, we estimate $f\left(x\right)$ and $f^{\left(2\right)}\left(x\right)$ by
\begin{eqnarray*}
\widetilde{f}_n\left(x\right)&=&\frac{1}{nb_n}\sum_{i=1}^n\delta_i\pi_{NWi}^{-1}K_b\left(\frac{x-X_i}{b_n}\right),\label{fn:norec:loc}\\
\widetilde{f}^{\left(2\right)}_n\left(x\right)&=&\frac{1}{nb_n^{\prime 3}}
\sum_{i=1}^n\delta_i\pi_{NWi}^{-1}K_{b^{\prime}}^{\left(2\right)}\left(\frac{x-X_i}{b^{\prime}_n}\right)\label{f:second:norec:loc}.
\end{eqnarray*}
where $K_b$ and $K_{b^{\prime}}$ are a kernels, $b_n$ and $b^{\prime}_n$ are respectively the associated bandwidth given in~(\ref{eq:h:initial}).\\
Then the plug-in estimator of the bandwidth $\left(h_n\right)$ using the nonrecursive estimator~(\ref{eq:Nad}), to estimate locally the density $f$ at a point $x$ is given by
\begin{eqnarray}\label{hoptim:rose:plug:in:loc}
\left(\left(\frac{\widetilde{f}_n\left(x\right)}{\left(\widetilde{f}_n^{\left(2\right)}\left(x\right)\right)^2}\right)^{1/5}\left\{\frac{R\left(K\right)}{\mu_2^2\left(K\right)}\right\}^{1/5}\widehat{\pi}_{NWn}^{-1/5}n^{-1/5}\right).
\end{eqnarray}
\section*{Acknowledgements}

We are grateful to two reviewers and the editor for their helpful comments, which have led to this substantially improved version of the paper.

\section*{}
\makeatletter
\renewcommand{\@biblabel}[1]{}
\makeatother


\begin{thebibliography}{99}

\bibitem[{Bojanic and Seneta(1973)}]{Boj73}
{Bojanic, R.} and {Seneta, E.} (1973).
\newblock {A unified theory of regularly varying sequences}.
\newblock \textit{Math. Z.}, {\bf 134}, 91--106.

\bibitem[{Chan et al.(2010)Chan, Lee and Peng}]{Cha10}
{Chan, N. H.} {Lee, T. C. M.} and {Peng, L.} (2010). \newblock {On nonparametric local inference for density estimation}.
\newblock \textit{Comput. Statist. Data Anal.}, {\bf 54},~509--515.

\bibitem[{Delaigle and Gijbels(2004)}]{Del04}
{Delaigle, A.} and {Gijbels, I.} (2004). 
\newblock {Practical bandwidth selection in deconvolution kernel density estimation}.
\newblock \textit{Comput. Statist. Data Anal.}, {\bf 45}, 249--267.

\bibitem[{Dempster et al.(1977)Dempster, Laird and Rubin}]{Dem77}
{Dempster, A. P.} {Laird, N. M.} and {Rubin, D. B.} (1977).
\newblock {Maximum likelihood from incomplete data via the EM algorithm}.
\newblock \textit{J. R. Statist. Soc. Ser. B Statist. Methodol.}, {\bf 39}, 1--22.


\bibitem[{Dubnicka(2009)}]{Dub09}
{Dubnicka, S. R.} (2009).
\newblock {Kernel density estimation with missing data and auxiliary variables}.
\newblock \textit{Aust. N. Z. J. Stat.}, {\bf 25}, 1175--1179.

\bibitem[{Duin(1976)}]{Dui76}
{Duin, R. P.} (1976).
\newblock {On the choice of smoothing parameters for parzen estimators of probablility density functions}.
\newblock \textit{IEEE Trans. Comput.}, {\bf 25}, 1175--1179.

\bibitem[{Dutta(2014)}]{Dut14}
{Dutta, S.} (2014). \newblock {Local smoothing using the bootstrap}. \newblock \textit{Comm. Statist. Simulation Comput.},  {\bf 43},~378--389. 

\bibitem[{Erdman and Emerson(2007)}]{Erd07}
{Erdman, C.} and {Emerson, J. W.} (2007).
\newblock {An R package for performing a bayesian analysis of change point problem}.
\newblock \textit{J. Statist. Software}, {\bf 23}, 1--13.

\bibitem[{Fan et al.(1995)}]{Fan95}
Fan, J., Heckmann, N.~E. and Wand, M.P. (1995). \newblock {Local polynomial kernel regression for generalized linear models and quasi likelihood functions}. \newblock \textit{J. Amer. Statist. Assoc.}, {\bf 90}, 141--150.

\bibitem[{Galambos and Seneta(1973)}]{Gal73}
{Galambos, J.} and {Seneta, E.} (1973).
\newblock {Regularly varying sequences}.
\newblock \textit{Proc. Amer. Math. Soc.}, {\bf 41}, 110--116.

\bibitem[{Hazelton(1999)}]{Haz99}
{Hazelton,~M.} (1999). \newblock {An optimal local bandwidth selector for kernel density estimation}. \newblock \textit{J. Statist. Plann. Inference}, {\bf 77},~37--50.

\bibitem[{Horvitz and Thompson(1952)}]{Hor52}
{Horvitz, D. G.} and {Thompson, D. J.} (1952).
\newblock {A generalization of sampling without replacement from a finite universes}.
\newblock \textit{J. Amer. Statist. Assoc.}, {\bf 47}, 663--685.

\bibitem[{Little and Rubin(2002)}]{Lit02} 
{Little, R.~J.~A} (2002).
\newblock {Statistical analysis with missing data}.
\newblock \textit{2nd. edn. Hoboken}, John Wiley \& Sons.


\bibitem[Marron (1988)]{Mar88}
{Marron, J. S.} (1988).
\newblock {Automatic smoothing parameter selction: a survey}.
\newblock \textit{empec.}, {\bf 13},~187--208.

\bibitem[{Mokkadem and Pelletier(2007)}]{Mok07}
{Mokkadem, A.} and {Pelletier, M.} (2007).
\newblock {A companion for the Kiefer-Wolfowitz-Blum stochastic approximation algorithm}.
\newblock \textit{Ann. Statist.}, {\bf 35}, 1749--1772.


\bibitem[{Mokkadem et~al.(2009)Mokkadem, Pelletier and Slaoui(2009)}]{Mok09} 
{Mokkadem, A.} {Pelletier, M.} and {Slaoui, Y.} (2009).
\newblock {The stochastic approximation method for the estimation of a multivariate probability density}.
\newblock \textit{J. Statist. Plann. Inference}, {\bf 139}, 2459--2478.

\bibitem[Nadaraya(1964)]{Nad64} 
{Nadaraya, E. A.} (1964).
\newblock {On estimating regression}.
\newblock \textit{Theory Probab. Appl.}, {\bf 10}, 186--190.

\bibitem[Parzen(1962)]{Par62} 
{Parzen, E.} (1962). 
\newblock {On estimation of a probability density and mode}. 
\newblock \textit{Ann. Math. Statist.}, {\bf 33},~1065--1076.

\bibitem[Qi et al.(2005)]{Qi05} 
Qi, L., Wang, C.~Y. and Prentice, R.~L. (2005). \newblock {Weighted estimators for proportional hazards regression with missing covariates}. \newblock \textit{J. Amer. Statist. Assoc.}, {\bf 100},~1250--1263.

\bibitem[R\'ev\'esz(1973)]{Rev73}
{R\'ev\'esz, P.} (1973).
\newblock {Robbins-Monro procedure in a Hilbert space and its
application in the theory of learning processes I}.
\newblock \textit{Studia Sci. Math. Hung.}, {\bf 8},~391--398.

\bibitem[R\'ev\'esz(1977)]{Rev77}
{R\'ev\'esz, P.} (1977).
\newblock {How to apply the method of stochastic approximation in the non-parametric estimation of a regression function}.
\newblock \textit{Math. Operationsforsch. Statist., Ser. Statistics.}, {\bf 8},~119--126.

\bibitem[Rosenblatt(1956)]{Ros56} 
{Rosenblatt, M.} (1956). 
\newblock {Remarks on some nonparametric estimates of a density function}. 
\newblock \textit{Ann. Math. Statist.}, {\bf 27},~832--837.

\bibitem[{Rosenbaum and Rubin(1983)}]{Ros83}
{Rosenbaum, P.R.} and {Rubin, D.B.} (1983). 
\newblock {The central role of the propensity score in observational studies
for causal effects}. 
\newblock \textit{Biometrika}, {\bf 70},~41--55.

\bibitem[Rudemo(1982)]{Rud82}
{Rudemo, M.} (1982).
\newblock {Empirical choice of histograms and kernel density estimators}.
\newblock \textit{Scand. J. Stat.}, {\bf 9},~65--78.

\bibitem[{Scott and Terrell(1987)}]{Sco87}
{Scott, D.W.} and {Terrell, G.R.} (1987). 
\newblock {Biased and unbiased cross-validation in density estimation}. 
\newblock \textit{J. Amer. Math. Soc.}, {\bf 82},~1131--1146.


\bibitem[{Sheather(1983)}]{She83}
{Sheather, S.~J.} (1983). \newblock {A data based algorithm for choosing the window when estimating the density at a point}. \newblock \textit{Comput. Statist. Data Anal.}, {\bf 1},~229--239.

\bibitem[{Sheather(1986)}]{She86}
{Sheather, S.~J.} (1986). \newblock {An improved data based algorithm for choosing the window when estimating the density at a point}. \newblock \textit{Comput. Statist. Data Anal.}, {\bf 4},~61--65.

\bibitem[{Silverman(1986)}]{Sil86} 
{Silverman, B.~W.} (1986).
\newblock {Density estimation for statistics and data analysis}.
\newblock \textit{Chapman and Hall}, London.

\bibitem[{Snijders et al.(2001)}]{Sni01} 
Snijders AM, Nowak N, Segraves R, Blackwood S, Brown N, Conroy J, Hamilton G, Hindle AK, Huey B, Kimura K, Law S, Myambo K, Palmer J, Ylstra B, Yue JP, Gray JW, Jain AN, Pinkel D, , Albertson DG (2001). 
\newblock {Assembly of Microarrays for Genome-wide Measurement
of DNA Copy Number}.
\newblock \textit{Nature Genetics}, {\bf 29},~263--246.

\bibitem[{Slaoui(2013)}]{Sla13} 
{Slaoui, Y.} (2013).
\newblock {Large and moderate principles for recursive kernel density estimators defined by stochastic approximation method}.
\newblock \textit{Serdica Math. J.}, {\bf 39},~53--82.

\bibitem[{Slaoui(2014a)}]{Sla14a} 
{Slaoui, Y.} (2014a).
\newblock {Bandwidth selection for recursive kernel density estimators defined by stochastic approximation method}.
\newblock \textit{J. Probab. Stat}, {\bf 2014}, ID 739640, doi:10.1155/2014/739640.

\bibitem[{Slaoui(2014b)}]{Sla14b} 
{Slaoui, Y.} (2014b).
\newblock {The stochastic approximation method for the estimation of a distribution function}.
\newblock \textit{Math. Methods Statist.}, {\bf 23},~306--325.


\bibitem[{Slaoui(2015)}]{Sla15} 
{Slaoui, Y.} (2015).
\newblock {Plug-In Bandwidth selector for recursive kernel regression estimators defined by stochastic approximation method.}
\newblock \textit{Stat. Neerl.}, {\bf 69},~483--509.

\bibitem[{Slaoui(2016)}]{Sla16} 
{Slaoui, Y.} (2016).
\newblock {Optimal bandwidth selection for semi-recursive kernel regression estimators.}
\newblock \textit{Stat. Interface}, {\bf 9},~375--388.


\bibitem[Thi\'ebault et al.(2000)]{Thi00}  
{Thi\'ebault, R., Morlat, P., Jacqmin-Gadda, H., Neau, D., Merci\'e, P., Dabis, F., and Ch\^ene, G. for the GECSA} (2000). 
\newblock {Clinical progression of HIV-1 infected patients according to the viral response during the first year of anti-retroviral treatment.} \newblock \textit{AIDS}, {\bf 14}:~971--978.


\bibitem[{Thombs and Sheather(1992)}]{Tho92} 
{Thombs, L. A., Sheather, S. J.} (1992). \newblock {Local bandwidth selection for density estimation. In: Proceedings
of the 22nd Symposium on the Interface.} \newblock \textit{New York: Springer},~111--116.

\bibitem[Tibshirani and Hastie(1987)]{Tib87} 
Tibshirani, R. and Hastie, T. (1987). \newblock {Local likelihood estimation}. \newblock \textit{J. Amer. Statist. Assoc.}, {\bf 82}, 559--567.

\bibitem[{Tsybakov(1990)}]{Tsy90} 
{Tsybakov, A.~B.} (1990).
\newblock {Recurrent estimation of the mode of a multidimensional distribution}.
\newblock \textit{Probl. Inf. Transm.}, {\bf 8},~119--126.


\bibitem[{Venkatraman and Olshen(2007)}]{Ven07}
{Venkatraman, E.S.} and {Olshen, A.B.} (2007). 
\newblock {A faster circular binary segmentation algorithm for the analysis of array CGH data}. 
\newblock \textit{Bioinformatics}, {\bf 23},~657--663.

\bibitem[Watson(1964)]{Wat64} 
{Watson, G. S.} (1964). 
\newblock {Smooth regression analysis}. 
\newblock \textit{Sankhya A.}, {\bf 26}, 359--372.
\end{thebibliography}
\end{document}